\documentclass{article}

\usepackage[titletoc,title]{appendix}
\usepackage[cp1250]{inputenc}
\usepackage[IL2]{fontenc}
\usepackage{yfonts,fancyhdr}
\usepackage{a4wide}
\usepackage[english]{babel}
\usepackage{euscript}
\usepackage{amstext,amsbsy,amscd,amssymb}
\usepackage{amsmath,enumerate}
\usepackage{amsfonts}
\usepackage{mathrsfs,tensor,xy}
\usepackage{graphicx}
\usepackage{microtype}
\usepackage{graphics}
\usepackage{tikz,caption,subcaption}
\usetikzlibrary{arrows,decorations.pathmorphing,backgrounds,positioning,fit,petri}
\include{diagxy}

\let\rarr=\rightarrow

\let\veps=\varepsilon
\let\mcal=\mathcal
\let\mfrak=\mathfrak
\let\eus=\EuScript

\def\N{\mathbb{N}}
\def\Z{\mathbb{Z}}

\def\C{\mathbb{C}}

\def\Mod{\mathop {\rm Mod} \nolimits}

\def\Hom{\mathop {\rm Hom} \nolimits}

\def\Aut{\mathop {\rm Aut} \nolimits}

\def\ad{\mathop {\rm ad} \nolimits}

\def\SL{\mathop {\rm SL} \nolimits}
\def\diag{\mathop {\rm diag} \nolimits}

\def\Ad{\mathop {\rm Ad} \nolimits}

\newsymbol\squares 1003

\long\def\proof #1{\noindent \emph{Proof.}\ #1 \hfill $\squares$
\medskip}

\newcounter{num}[section]
\numberwithin{equation}{section}
\numberwithin{num}{section}

\long\def\definition #1 {\refstepcounter{num} \noindent {\bf
Definition \thenum.} #1

\medskip}

\long\def\lemma #1{\refstepcounter{num}  \noindent {\bf Lemma
\thenum.} #1

\medskip}

\newenvironment{enum}{\begin{list}{}{\topsep=2pt \itemsep=0pt
\parsep=0pt}}{\end{list}}

\newcommand*\riso{%
  \xrightarrow[]{\raisebox{-0.25em}{\smash{\ensuremath{\sim}}}}%
}

\makeatletter
\newcommand*\if@single[3]{%
  \setbox0\hbox{${\mathaccent"0362{#1}}^H$}%
  \setbox2\hbox{${\mathaccent"0362{\kern0pt#1}}^H$}%
  \ifdim\ht0=\ht2 #3\else #2\fi
  }
\newcommand*\rel@kern[1]{\kern#1\dimexpr\macc@kerna}
\newcommand*\widebar[1]{\@ifnextchar^{{\wide@bar{#1}{0}}}{\wide@bar{#1}{1}}}
\newcommand*\wide@bar[2]{\if@single{#1}{\wide@bar@{#1}{#2}{1}}{\wide@bar@{#1}{#2}{2}}}
\newcommand*\wide@bar@[3]{%
  \begingroup
  \def\mathaccent##1##2{%
    \if#32 \let\macc@nucleus\first@char \fi
    \setbox\z@\hbox{$\macc@style{\macc@nucleus}_{}$}%
    \setbox\tw@\hbox{$\macc@style{\macc@nucleus}{}_{}$}%
    \dimen@\wd\tw@
    \advance\dimen@-\wd\z@
    \divide\dimen@ 3
    \@tempdima\wd\tw@
    \advance\@tempdima-\scriptspace
    \divide\@tempdima 10
    \advance\dimen@-\@tempdima
    \ifdim\dimen@>\z@ \dimen@0pt\fi
    \rel@kern{0.6}\kern-\dimen@
    \if#31
      \overline{\rel@kern{-0.6}\kern\dimen@\macc@nucleus\rel@kern{0.4}\kern\dimen@}%
      \advance\dimen@0.4\dimexpr\macc@kerna
      \let\final@kern#2%
      \ifdim\dimen@<\z@ \let\final@kern1\fi
      \if\final@kern1 \kern-\dimen@\fi
    \else
      \overline{\rel@kern{-0.6}\kern\dimen@#1}%
    \fi
  }%
  \macc@depth\@ne
  \let\math@bgroup\@empty \let\math@egroup\macc@set@skewchar
  \mathsurround\z@ \frozen@everymath{\mathgroup\macc@group\relax}%
  \macc@set@skewchar\relax
  \let\mathaccentV\macc@nested@a
  \if#31
    \macc@nested@a\relax111{#1}%
  \else
    \def\gobble@till@marker##1\endmarker{}%
    \futurelet\first@char\gobble@till@marker#1\endmarker
    \ifcat\noexpand\first@char A\else
      \def\first@char{}%
    \fi
    \macc@nested@a\relax111{\first@char}%
  \fi
  \endgroup
}
\makeatother

\title{On the composition structure of the twisted Verma modules for $\mathfrak{sl}(3,\mathbb{C})$}

\author{Libor Křižka, Petr Somberg}

\AtEndDocument{\bigskip{\footnotesize%
  (L.\,Křižka) \textsc{Mathematical Institute of Charles University, Sokolovská 83, 180\,00 Praha 8, Czech Republic} \par
  \textit{E-mail address}: \texttt{krizka@karlin.mff.cuni.cz} \par
  \addvspace{\medskipamount}
  (P.\,Somberg) \textsc{Mathematical Institute of Charles University, Sokolovská 83, 180\,00 Praha 8, Czech Republic} \par
  \textit{E-mail address}: \texttt{somberg@karlin.mff.cuni.cz}
}}

\begin{document}
\date{}
\maketitle

\begin{abstract}

We discuss some aspects of the composition structure of twisted Verma modules
for the Lie algebra $\mathfrak{sl}(3, \mathbb{C})$, including the explicit
structure of singular vectors for both $\mathfrak{sl}(3, \mathbb{C})$ and one
of its Lie subalgebras $\mathfrak{sl}(2, \mathbb{C})$, and also of their generators.
Our analysis is based on the use of partial Fourier tranform applied to the
realization of twisted Verma modules as ${\fam2 D}$-modules on the Schubert cells
in the full flag manifold for $\SL(3, \mathbb{C})$.


\end{abstract}

\thispagestyle{empty}

\tableofcontents

\section*{Introduction}
\addcontentsline{toc}{section}{Introduction}

The objects of central interest in the representation theory
of complex simple Lie algebras are the Harish-Chandra modules.
It is well known that there is a categorical equivalence between
principal series Harish-Chandra modules and twisted Verma modules
as objects of the Bernstein-Gelfand-Gelfand category $\mcal{O}$.
The twisted Verma modules are studied from various perspectives
including the Lie algebra (co)homology of the twisted nilradical,
the Schubert cell decomposition of full flag manifolds and algebraic
techniques of twisting functors applied to Verma modules, in
\cite{Feigin-Frenkel1990}, \cite{Andersen-Lauritzen2003}, \cite{Soergel1998}
and references therein.

Combinatorial conditions for the existence of homomorphisms between twisted
Verma modules were studied in \cite{Abe2011},
but there is basically no information on precise positions and properties
of elements responsible for a non-trivial composition structure of twisted
Verma modules. The modest aim of the present article
is the study of some aspects related to the composition structure
of twisted Verma modules for the Lie algebra $\mathfrak{sl}(3, \C)$
by geometrical methods,
through their realization as $\mcal{D}$-modules supported on Schubert cells,
cf.\ \cite{Hotta-book}, \cite{Beilinson-Bernstein1981}.
Namely, we discuss in the case of $\mathfrak{sl}(3, \C)$ a few results parallel to
the development for (untwisted) generalized Verma modules in
\cite{Kobayashi2012}, \cite{koss}.

Let us briefly describe the content of our article.
First of all, in Section \ref{sec:intro} we briefly review various characterizing
properties of twisted Verma modules compared to the
untwisted Verma modules. Based on the
action of a simple Lie algebra on its full flag manifold, see
e.g. \cite{Hotta-book}, \cite{Krizka-Somberg2015} for rather explicit
description, in Section \ref{sec:twisted}
we write down
explicit realizations of the highest weight twisted
$\mathfrak{sl}(3,\C)$-Verma modules. For all twistings realized
by elements $w$ of the Weyl group $W$, the isomorphism given by
a partial Fourier transform allows us to analyze several basic questions
on twisted Verma modules not accessible in the literature.
We shall carry out this procedure for the structure of singular vectors
and the generators of $\mfrak{sl}(3,\C)$-Verma modules twisted by
$w=s_1$. Another our result concerns the application of ideas on
the decomposition of twisted Verma modules with respect to a
reductive Lie subalgebra $\mathfrak{sl}(3,\C)$, thereby generalizing
the results analogous to \cite{koss} towards the twisted Verma modules.
Here we consider the simplest example of an embedded
$\mathfrak{sl}(2,\C)\subset\mathfrak{sl}(3,\C)$ and
produce a complete list of singular vectors responsible
for the branching problem of twisted $\mathfrak{sl}(3,\C)$-Verma module.
In our situation we also observe that the $s_1$-twisted $\mfrak{sl}(3,\C)$-Verma modules
are generated by single vector (which is not of highest weight),
a property analogous to the case of (untwisted) Verma modules.
In the last Section \ref{out} we highlight our results in the framework
of (un)known properties of the objects of the Bernstein-Gelfand-Gelfand
category $\mcal{O}$ (see e.g.\ \cite{Humphreys-book}).

\section{Twisted Verma modules and their characterizations}
\label{sec:intro}

Let $G$ be a connected complex semisimple Lie group, $H \subset G$
a maximal torus of $G$, $B \subset G$ a Borel subgroup of $G$ containing
$H$, and $W=N_G(H)/H$ the Weyl group of $G$. Furthermore, let
$\mfrak{g}$, $\mfrak{h}$ and $\mfrak{b}$ be the Lie algebras of $G$,
$H$ and $B$, respectively. Finally, let $\mathfrak{n}$ be the positive
nilradical of the Borel subalgebra $\mathfrak{b}$ and $\widebar{\mathfrak{n}}$
the opposite (negative) nilradical. We denote by $N$ and $\widebar{N}$ the Lie subgroups
 of $G$ corresponding to the Lie subalgebras $\mathfrak{n}$ and
$\widebar{\mathfrak{n}}$, respectively.

The objects of our interest are the twisted Verma modules
$M^w_\mathfrak{g}(\lambda)$, parametrized by $\lambda\in\mathfrak{h}^*$
and the twisting $w\in W$. The twisted Verma modules $M^w_\mathfrak{g}(\lambda)$ have for all $w\in W$ the same
character as the Verma module $M^\mfrak{g}_\mfrak{b}(\lambda)$ induced from the $1$-dimensional
$\mathfrak{b}$-module $\C_\lambda$,
\begin{align*}
M^\mathfrak{g}_\mathfrak{b}(\lambda)\equiv M^e_\mathfrak{g}(\lambda)=
U(\mathfrak{g})\otimes_{U(\mathfrak{b})}\!\C_\lambda,
\end{align*}
with highest weight $\lambda\in\mathfrak{h}^*$ and $e\in W$. However, the extensions of
simple sub-quotients in twisted
Verma modules differ from extensions in Verma modules. As $U(\mathfrak{g})$-modules
they are objects of the Bernstein-Gelfand-Gelfand category $\mcal{O}$,
i.e.\ finitely generated $U(\mfrak{g})$-modules, $\mfrak{h}$-semisimple and locally $\mfrak{n}$-finite.

Let us denote by $e$ and $w_0$ the identity and the longest element of $W$,
respectively, and let $\ell \colon W \rarr  \N_0$ be the length function on $W$.
The Weyl group $W$ acts by $\rho$-affine action on $\mfrak{h}^*$,
$w\cdot \lambda=w(\lambda+\rho)-\rho$, and gives four Lie subalgebras
\begin{enumerate}
\item[1)]
$\widebar{\mfrak{n}}=\mfrak{n}_w^-\oplus \widebar{\mfrak{n}}_w^-$:
\quad
$\mfrak{n}_w^-=\widebar{\mfrak{n}} \cap \Ad(\dot{w})(\mfrak{n}),\
\widebar{\mfrak{n}}_w^-=\widebar{\mfrak{n}} \cap \Ad(\dot{w})(\widebar{\mfrak{n}})$,
\item[2)]
${\mfrak{n}}=\mfrak{n}_w^+\oplus \widebar{\mfrak{n}}_w^+$: \quad
$\mfrak{n}_w^+={\mfrak{n}} \cap \Ad(\dot{w})(\mfrak{n}),\
\widebar{\mfrak{n}}_w^+={\mfrak{n}} \cap \Ad(\dot{w})(\widebar{\mfrak{n}})$.
\end{enumerate}
The universal enveloping algebra
$U(\mfrak{n}_w^-)$ is a graded subalgebra of $U(\widebar{\mfrak{n}})$,
determined by $U(\mfrak{n}_w^-)_0=\C$, $U(\mfrak{n}_w^-)_{-1}=\mfrak{n}_w^-$
for all $w\in W$, and $U(\mfrak{n}_e^-)=\C$,
$U(\mfrak{n}_{w_0}^-)=U(\widebar{\mfrak{n}})$. The graded dual of $U(\mfrak{n}_w^-)$
is defined by
$(U(\mfrak{n}_w^-))^*_n=\Hom_\C((U(\mfrak{n}_w^-))_{-n}, \C)$
for all $n\in\Z$.
\smallskip

There are several equivalent characterizing properties of twisted Verma modules $M^w_\mathfrak{g}(\lambda)$
for $\lambda \in \mfrak{h}^*$ and $w \in W$, see \cite{Andersen-Lauritzen2003}, \cite[Chapter 12]{Humphreys-book}
for detailed discussion.
\begin{enumerate}
\item[1)]
The Lie algebra cohomology of the twisted opposite nilradical
$\widebar{\mathfrak{n}}_w^+\oplus\widebar{\mathfrak{n}}_w^- =\Ad(\dot{w})(\widebar{\mathfrak{n}})$
with coefficients in $M^w_\mathfrak{g}(\lambda)$
is
\begin{align}
H^i(\widebar{\mathfrak{n}}_w^+\oplus\widebar{\mathfrak{n}}_w^-, M^w_\mathfrak{g}(\lambda))\simeq
\begin{cases}
  \C_{\lambda +w(\rho)+\rho} & \text{if $i=\dim \mfrak{n} - \ell(w)$}, \\
  0 & \text{if $i \neq \dim \mfrak{n} - \ell(w)$}
\end{cases}
\end{align}
as $\mathfrak{h}$-modules. In particular, $M^w_\mathfrak{g}(\lambda)$ is
a free $U(\Ad(\dot{w})(\widebar{\mathfrak{n}}) \cap \widebar{\mfrak{n}})$-module,
while its graded dual $(M^w_\mathfrak{g}(\lambda))^*$ is a free
$U(\Ad(\dot{w})(\widebar{\mathfrak{n}}) \cap \mfrak{n})$-module.
\item[2)]
Let us consider the full flag manifold $G/B$ and the Schubert cell
$X_w$ for $w\in W$ defined as the $N$-orbit $X_w=NwB/B\subset G/B$,
where $\dim X_w=\ell(w)$. Then there is an isomorphism of
$U(\mathfrak{g})$-modules for the local, relative to $X_w$,
sheaf cohomology of a homogeneous vector bundle $\mcal{L}(\lambda)$,
\begin{align}
H^i_{X_w}(G/B,\mcal{L}(\lambda))\simeq \begin{cases}
  M^w_\mfrak{g}(ww_0 \cdot \lambda)  & \text{if $i=\dim \mfrak{n} - \ell(w)$}, \\
  0 & \text{if $i \neq \dim \mfrak{n} - \ell(w)$.}
\end{cases}
\end{align}
In particular, the Verma modules are supported on the closed Schubert cell while the
contragradient Verma modules are supported on the open (dense)
Schubert cell.
\item[3)]
For $w\in W$, the $U(\mfrak{g})$-bimodule $S_w=U(\mfrak{g})\otimes_{U(\mfrak{n}_w^-)}\!(U(\mfrak{n}_w^-))^*$
allows to define a functor $T_w \colon  \mcal{O} \rarr \mcal{O}$ (called twisting functor)
by
\begin{align}
T_w \colon  M\mapsto \varphi_w(S_w\otimes_{U(\mathfrak{g})}M).
\end{align}
Here $\varphi_w =\Ad(\dot{w}^{-1}) \colon \mathfrak{g} \rarr \Aut(\mathfrak{g})$
indicates the conjugation of the action by $\mathfrak{g}$ on
the twisted module. In particular, we have
$M^w_\mathfrak{g}(\lambda)=T_w(M^\mfrak{g}_\mfrak{b}(w\cdot \lambda))$.
\end{enumerate}

Twisted Verma modules $\smash{M^w_\mfrak{g}}(\lambda)$ for $\lambda \in \mfrak{h}^*$ and $w \in W$
can be realized in the framework of $\mcal{D}$-modules on the flag manifold $X=G/B$.
There is a $G$-equivariant sheaf of rings of twisted differential operators
$\mcal{D}^\lambda_X$ on $X$, see \cite{Beilinson-Bernstein1981}, \cite{Kashiwara1989}, which is
for an integral dominant weight $\lambda+\rho$ a sheaf of rings of differential operators acting on $\mcal{L}(\lambda+\rho)$.
The $G$-equivariance of $\mcal{D}^\lambda_X$ ensures
the existence of a Lie algebra morphism
\begin{align}
  \alpha_\lambda \colon \mfrak{g} \rarr \Gamma(X,\mcal{D}^\lambda_X)
\end{align}
and a localization functor
\begin{align}
  \Delta \colon \Mod(\mfrak{g}) \rarr \Mod(\mcal{D}^\lambda_X).
\end{align}
Then the $\mcal{D}^\lambda_X$-module $\Delta(\smash{M^w_\mfrak{g}}(w\cdot(\lambda-\rho)))$ is realized in
the vector space of distributions supported on the Schubert cell $X_w$ of $X$, see \cite[Chapter 11]{Frenkel-Ben-Zvi-book}.

\section{Twisted Verma modules for $\mfrak{sl}(3,\C)$}
\label{sec:twisted}

We shall consider the complex semisimple Lie group $G=\SL(3,\C)$ given by $3\times 3$
complex matrices of unit determinant and its Lie algebra $\mfrak{g}=\mfrak{sl}(3,\C)$.
The Cartan subalgebra $\mfrak{h}$ of $\mfrak{g}$ is given by diagonal matrices
$\mfrak{h}=\{\diag(a_1,a_2,a_3);\, a_1,a_2,a_3\in \C,\ a_1+a_2+a_3=0\}$.
For $i=1,2,3$, we define $\veps_i \in \mfrak{h}^*$ by $\veps_i(\diag(a_1,a_2,a_3))=a_i$.
Then the root system of $\mfrak{g}$ with respect to $\mfrak{h}$ is
$\Delta =\{\veps_i-\veps_j;\, 1\leq i\neq j\leq 3\}$, the positive root system is
$\Delta^+=\{\veps_i-\veps_j;\, 1\leq i < j\leq 3\}$ and the set of simple roots is
$\Pi=\{\alpha_1,\alpha_2\}$, $\alpha_1=\veps_1-\veps_2$, $\alpha_2=\veps_2-\veps_3$.
The fundamental weights are $\omega_1=\veps_1$, $\omega_2=\veps_1+\veps_2$, and the
smallest regular integral dominant weight is $\rho=\omega_1+\omega_2$. The notation
$\lambda=(\lambda_1, \lambda_2)$ means $\lambda=\lambda_1\omega_1+\lambda_2\omega_2$.

We choose the basis of root spaces of $\mfrak{g}$ as
\begin{align*}
  f_1&=f_{\alpha_1}=
  \begin{pmatrix}
    0 & 0 & 0 \\
    1 & 0 & 0 \\
    0 & 0 & 0
  \end{pmatrix}\!, &
	f_2&=f_{\alpha_2}=
  \begin{pmatrix}
    0 & 0 & 0  \\
    0 & 0 & 0 \\
    0 & 1 & 0
  \end{pmatrix}\!, &
	f_{12}&=f_{\alpha_1+\alpha_2}=
  \begin{pmatrix}
    0 & 0 & 0  \\
    0 & 0 & 0 \\
    1 & 0 & 0
  \end{pmatrix}\!,	\\
   e_1&=e_{\alpha_1}=
  \begin{pmatrix}
    0 & 1 & 0 \\
    0 & 0 & 0 \\
    0 & 0 & 0
  \end{pmatrix}\!, &
  	e_2&=e_{\alpha_2}=
  \begin{pmatrix}
    0 & 0 & 0  \\
    0 & 0 & 1 \\
    0 & 0 & 0
  \end{pmatrix}\!, &
	e_{12}&=e_{\alpha_1+\alpha_2}=
  \begin{pmatrix}
    0 & 0 & 1  \\
    0 & 0 & 0 \\
    0 & 0 & 0
  \end{pmatrix}\!,
\end{align*}
and the basis of the Cartan subalgebra $\mfrak{h}$ is given by coroots
\begin{align*}
	h_1&=h_{\alpha_1}=
  \begin{pmatrix}
    1 & 0 & 0 \\
    0 & -1 & 0 \\
    0 & 0 & 0
  \end{pmatrix}\!, &
	h_2&=h_{\alpha_2}=
  \begin{pmatrix}
    0 & 0 & 0  \\
    0 & 1 & 0 \\
    0 & 0 & -1
  \end{pmatrix}\!.
\end{align*}
These matrices fulfill, among others, the commutation relations
$[f_{\alpha_1}, f_{\alpha_2}]=-f_{\alpha_1+\alpha_2}$
and $[e_{\alpha_1}, e_{\alpha_2}]=e_{\alpha_1+\alpha_2}$.

The Weyl group $W$ of $G$ is generated by simple reflections
$s_1=s_{\alpha_1}$ and $s_2=s_{\alpha_2}$, where the action of
$W$ on $\mfrak{h}^*$ is given by
\begin{align*}
s_1(\alpha_1)=-\alpha_1,\quad s_1(\alpha_2)=\alpha_1+\alpha_2,\quad
s_2(\alpha_1)=\alpha_1+\alpha_2,\quad s_2(\alpha_2)=-\alpha_2,
\end{align*}
and $|W|=6$ with $W=\{e, s_1, s_2, s_1s_2, s_2s_1, s_1s_2s_1=s_2s_1s_2 \}$.
Consequently, there are six Schubert cells $X_w$ in $G/B$ isomorphic to $X_w \simeq \C^{\ell(w)}$:
\begin{align*}
 \dim(X_e)=0,\ \dim(X_{s_1})=\dim(X_{s_2})=1,\ \dim(X_{s_1s_2})=\dim(X_{s_2s_1})=2,\
\dim(X_{s_1s_2s_1})=3.
\end{align*}
For the representatives of the elements of $W$ in $G$ we take the matrices
\begin{align*}
  \dot{e}&=
  \begin{pmatrix}
    1 & 0 & 0 \\
    0 & 1 & 0 \\
    0 & 0 & 1
  \end{pmatrix}\!, &
  \dot{s}_1&=
  \begin{pmatrix}
    0 & 1 & 0  \\
    -1 & 0 & 0 \\
    0 & 0 & 1
  \end{pmatrix}\!, &
	\dot{s}_2&=
  \begin{pmatrix}
    1 & 0 & 0  \\
    0 & 0 & 1 \\
    0 & -1 & 0
  \end{pmatrix}\!,	 \\
 \dot{s_1}\dot{s_2}&=
  \begin{pmatrix}
    0 & 0 & 1 \\
    -1 & 0 & 0 \\
    0 & -1 & 0
  \end{pmatrix}\!, &
	\dot{s}_2\dot{s}_1&=
  \begin{pmatrix}
    0 & 1 & 0  \\
    0 & 0 & 1 \\
    1 & 0 & 0
  \end{pmatrix}\!, &
	\dot{s}_1\dot{s}_2\dot{s}_1&=
  \begin{pmatrix}
    0 & 0 & 1  \\
    0 & -1 & 0 \\
    1 & 0 & 0
  \end{pmatrix}\!.
\end{align*}


We denote by $(x,y,z)$ the linear coordinate functions on $\widebar{\mfrak{n}}$
with respect to the basis $(f_1, f_2, f_{12})$ of the opposite nilradical
$\widebar{\mfrak{n}}$, and by $(\xi_x,\xi_y,\xi_z)$ the dual linear coordinate
functions on $\widebar{\mfrak{n}}^*$.

Let us consider the partial dual space $\widebar{\mfrak{n}}^{*,w}$
of $\widebar{\mfrak{n}}$ defined by
\begin{align}
  \widebar{\mfrak{n}}^{*,w} = (\widebar{\mfrak{n}}_{\smash{w^{-1}}}^-)^* \oplus \mfrak{n}_{\smash{w^{-1}}}^-\, ,
\end{align}
so that
\begin{align}
  (\xi_{x_\alpha},\,\alpha\in w^{-1}(\Delta^+)\cap\Delta^+,\ x_\alpha,\,\alpha\in w^{-1}(-\Delta^+)\cap\Delta^+)
\end{align}
with $x_{\alpha_1}=x$, $x_{\alpha_2}=y$, $x_{\alpha_1+\alpha_2}=z$ are linear coordinate functions on $\widebar{\mfrak{n}}^{*,w}$.
Moreover, the Weyl algebra $\eus{A}^\mfrak{g}_{\widebar{\mfrak{n}}}$
of $\widebar{\mfrak{n}}$ is generated by $\{x,\partial_x,y,\partial_y,z,\partial_z\}$,
and the Weyl algebra $\eus{A}^\mfrak{g}_{\widebar{\mfrak{n}}^{*,w}}$ of $\widebar{\mfrak{n}}^{*,w}$  is generated by
\begin{align}
\{\xi_{x_\alpha},\partial_{\xi_{x_\alpha}},\, \alpha\in w^{-1}(\Delta^+)\cap\Delta^+,\ x_\alpha,\partial_{x_\alpha},\, \alpha\in w^{-1}(-\Delta^+)\cap\Delta^+\}.
\label{eq:weyl algebra dual gen}
\end{align}
There is a canonical isomorphism
$\mcal{F}^w \colon \eus{A}^\mfrak{g}_{\widebar{\mfrak{n}}} \rarr \eus{A}^\mfrak{g}_{\widebar{\mfrak{n}}^{*,w}}$
of associative $\C$-algebras called the partial Fourier transform,
defined with respect to the generators \eqref{eq:weyl algebra dual gen} by
\begin{align}
\begin{aligned}
x_\alpha &\mapsto -\partial_{\xi_{x_\alpha}},&   \partial_{x_\alpha} \mapsto \xi_{x_\alpha},
 & \quad \text{for}\ \alpha\in w^{-1}(\Delta^+)\cap\Delta^+, \\
x_\alpha &\mapsto x_\alpha, & \partial_{x_\alpha} \mapsto \partial_{x_\alpha},
 & \quad \text{for}\ \alpha\in w^{-1}(-\Delta^+)\cap\Delta^+.
\end{aligned}
\end{align}
The partial Fourier transform is independent of the choice of linear coordinates on $\widebar{\mfrak{n}}$.

The Verma modules $M^\mfrak{g}_\mfrak{b}(\lambda-\rho)$, $\lambda \in \mfrak{h}^*$,
can be realized as $\eus{A}^\mfrak{g}_{\widebar{\mfrak{n}}}/I_e$ for
$I_e$ the left ideal of $\eus{A}^\mfrak{g}_{\widebar{\mfrak{n}}}$ defined by $I_e=(x,y,z)$,
see e.g.\ \cite{Krizka-Somberg2015}. The structure
of $\mfrak{g}$-module on $\eus{A}^\mfrak{g}_{\widebar{\mfrak{n}}}/I_e$ is realized through the embedding
$\pi_\lambda \colon \mfrak{g} \rarr \eus{A}^\mfrak{g}_{\widebar{\mfrak{n}}}$ given by
\begin{align}
\pi_\lambda(X)= -\sum_{\alpha \in \Delta^+}\bigg[{\ad(u(x))e^{\ad(u(x))} \over e^{\ad(u(x))}-{\rm id}_{\widebar{\mfrak{n}}}}\,(e^{-\ad(u(x))}X)_{\widebar{\mfrak{n}}}\bigg]_\alpha \partial_{x_\alpha} + (\lambda+\rho)((e^{-\ad(u(x))}X)_\mfrak{b})
\end{align}
for all $X\in \mfrak{g}$, where $[Y]_\alpha$ denotes the $\alpha$-th coordinate
of $Y \in \widebar{\mfrak{n}}$ with respect to the basis $(f_\alpha;\, \alpha \in \Delta^+)$
of $\widebar{\mfrak{n}}$ and $u(x)=\sum_{\alpha \in \Delta^+} x_\alpha f_\alpha$.
The twisted Verma modules are realized by $\eus{A}^\mfrak{g}_{\widebar{\mfrak{n}}}/I_w$, where
$I_w$ is the left ideal of $\eus{A}^\mfrak{g}_{\widebar{\mfrak{n}}}$ defined by
\begin{align}
I_w = (x_\alpha,\,\alpha\in w^{-1}(\Delta^+)\cap\Delta^+,\ \partial_{x_\alpha},\,\alpha\in w^{-1}(-\Delta^+)\cap\Delta^+)
\end{align}
with $x_{\alpha_1}=x$, $x_{\alpha_2}=y$, $x_{\alpha_1+\alpha_2}=z$.
The list of all possibilities looks as follows:
\begin{enumerate}
  \item[1)] $w=e$, $I_e=(x,y,z)$, $\eus{A}^\mfrak{g}_{\widebar{\mfrak{n}}}/I_e \simeq
\C[\partial_x, \partial_y, \partial_z]$;
  \item[2)] $w=s_1$,  $I_{s_1}=(\partial_x,y,z)$, $\eus{A}^\mfrak{g}_{\widebar{\mfrak{n}}}/I_{s_1}
\simeq \C[x, \partial_y, \partial_z]$;
  \item[3)] $w=s_2$, $I_{s_2}=(x,\partial_y,z)$, $\eus{A}^\mfrak{g}_{\widebar{\mfrak{n}}}/I_{s_2}
\simeq \C[\partial_x, y, \partial_z]$;
  \item[4)] $w=s_1s_2$, $I_{s_1s_2}=(x,\partial_y,\partial_z)$, $\eus{A}^\mfrak{g}_{\widebar{\mfrak{n}}}/I_{s_1s_2}
\simeq \C[\partial_x, y, z]$;
  \item[5)] $w=s_2s_1$, $I_{s_2s_1}=(\partial_x,y,\partial_z)$, $\eus{A}^\mfrak{g}_{\widebar{\mfrak{n}}}/I_{s_2s_1}
\simeq \C[x,\partial_y, z]$;
  \item[6)] $w=s_1s_2s_1$, $I_{s_1s_2s_1}=(\partial_x,\partial_y,\partial_z)$, $\eus{A}^\mfrak{g}_{\widebar{\mfrak{n}}}/I_{s_1s_2s_1}
\simeq \C[x, y, z]$.
\end{enumerate}
In particular, the twisted Verma modules are realized as
$M^w_\mfrak{g}(\lambda)\simeq \eus{A}^\mfrak{g}_{\widebar{\mfrak{n}}}/I_w$,
where $\pi^w_\lambda \colon \mfrak{g}\rarr \eus{A}^\mfrak{g}_{\widebar{\mfrak{n}}}/I_w$
is defined by
\begin{align}
\pi^w_\lambda=\pi_{w^{-1}(\lambda+\rho)}\circ \Ad(\dot{w}^{-1})
\end{align}
with $w^{-1}$ acting in the standard and not the $\rho$-shifted manner.
Since
$\mcal{F}^w \colon \eus{A}^\mfrak{g}_{\widebar{\mfrak{n}}} \rarr \eus{A}^\mfrak{g}_{\widebar{\mfrak{n}}^{*,w}}$
is an isomorphism of associative $\C$-algebras, the composition
\begin{align}
  \hat{\pi}_\lambda^w = \mcal{F}^w \circ \pi_\lambda^w
\end{align}
gives the homomorphism
$\hat{\pi}^w_\lambda \colon U(\mfrak{g}) \rarr \eus{A}^\mfrak{g}_{\widebar{\mfrak{n}}^{*,w}}$
of associative $\C$-algebras and the twisted Verma modules are realized as
$M^w_\mfrak{g}(\lambda) \simeq \eus{A}^\mfrak{g}_{\widebar{\mfrak{n}}^{*,w}}/\mcal{F}^w(I_w)$.

\subsection{(Untwisted) Verma modules}

Let us first consider the case of Verma modules.
The untwisted Verma module $M^\mfrak{g}_\mfrak{b}(\lambda)\equiv M_\mfrak{g}^{e}(\lambda)$ for $\lambda=\lambda_1\omega_1+\lambda_2\omega_2$ is isomorphic to
\begin{align}
M_{\mathfrak g}^{e}(\lambda)\simeq \eus{A}^\mfrak{g}_{\widebar{\mfrak{n}}}/I_e \simeq
\C[\partial_x, \partial_y, \partial_z], \quad I_e=(x, y, z),
\end{align}
where the embedding
$\pi^e_\lambda \colon  \mfrak{g}\rarr \eus{A}^\mfrak{g}_{\widebar{\mfrak{n}}}$
and so the $\mfrak{g}$-module structure on $\eus{A}^\mfrak{g}_{\widebar{\mfrak{n}}}/I_e$ are given by
\begin{align*}
\begin{aligned}
\pi^e_\lambda(f_1)&=-\partial_x+{\textstyle {1 \over 2}}y\partial_z, \\
\pi^e_\lambda(f_2)&=-\partial_y-{\textstyle {1 \over 2}}x\partial_z, \\
\pi^e_\lambda(f_{12})&=-\partial_z, \\ 
\pi^e_\lambda(e_1)&=x^2\partial_x + (z-{\textstyle {1\over 2}}xy)\partial_y+({\textstyle {1\over 4}}x^2y+ {\textstyle {1\over 2}}xz)\partial_z + (\lambda_1+2)x, \\
\pi^e_\lambda(e_2)&=y^2\partial_y -(z+{\textstyle {1\over 2}}xy)\partial_x -({\textstyle {1\over 4}}xy^2-{\textstyle {1\over 2}}yz)\partial_z + (\lambda_2+2)y, \\
\pi^e_\lambda(e_{12})&=(xz+{\textstyle {1\over 2}}x^2y)\partial_x + (yz-{\textstyle {1\over 2}}xy^2)\partial_y + (z^2+{\textstyle {1 \over 4}}x^2y^2)\partial_z + (\lambda_1+\lambda_2+4)z + {\textstyle {1\over 2}}(\lambda_1-\lambda_2)xy, \\
\pi^e_\lambda(h_1)&=2x\partial_x-y\partial_y+z\partial_z+\lambda_1+2, \\
\pi^e_\lambda(h_2)&=-x\partial_x+2y\partial_y+z\partial_z+\lambda_2+2. \label{nontwistedaction}
\end{aligned}
\end{align*}
This representation corresponds to the Verma module with the highest
weight $\lambda=(\lambda_1, \lambda_2)$, and the Fourier dual representation
acts in the Fourier dual variables on $\C[\xi_x, \xi_y, \xi_z]$ by
\begin{align*}
\begin{aligned}
\hat{\pi}^e_\lambda(f_1)&=-\xi_x-{\textstyle {1 \over 2}}\xi_z\partial_{\xi_y}, \\
\hat{\pi}^e_\lambda(f_2)&=-\xi_y+{\textstyle {1 \over 2}}\xi_z\partial_{\xi_x}, \\
\hat{\pi}^e_\lambda(f_{12})&=-\xi_z, \\
\hat{\pi}^e_\lambda(e_1)&=-\xi_y\partial_{\xi_z}
+(\xi_x\partial_{\xi_x}+{\textstyle {1 \over 2}}\xi_z\partial_{\xi_z}-\lambda_1)\partial_{\xi_x}
-{\textstyle {1 \over 2}}(\xi_{y}+{\textstyle {1 \over 2}}\xi_z\partial_{\xi_x})\partial_{\xi_x}\partial_{\xi_y}, \\
\hat{\pi}^e_\lambda(e_2)&=\xi_x\partial_{\xi_z}
+(\xi_y\partial_{\xi_y}+{\textstyle {1 \over 2}}\xi_z\partial_{\xi_z}-\lambda_2)\partial_{\xi_y}
-{\textstyle {1 \over 2}}(\xi_x-{\textstyle {1 \over 2}}\xi_{z}\partial_{\xi_y})\partial_{\xi_x}\partial_{\xi_y}, \\
\hat{\pi}^e_\lambda(e_{12})&=(\xi_x\partial_{\xi_x}+\xi_y\partial_{\xi_y}+\xi_z\partial_{\xi_z}
-\lambda_1-\lambda_2)\partial_{\xi_z}
-{\textstyle {1 \over 2}}(\xi_x\partial_{\xi_x}-\xi_y\partial_{\xi_y}
-\lambda_1+\lambda_2-{\textstyle {1 \over 2}}\xi_z\partial_{\xi_x}\partial_{\xi_y})\partial_{\xi_x}\partial_{\xi_y}, \\ \hat{\pi}^e_\lambda(h_1)&=-2\xi_x\partial_{\xi_x}+\xi_y\partial_{\xi_y}-\xi_z\partial_{\xi_z}+\lambda_1,
\\ \hat{\pi}^e_\lambda(h_2)&=\xi_x\partial_{\xi_x}-2\xi_y\partial_{\xi_y}-\xi_z\partial_{\xi_z}+\lambda_2.
\end{aligned}
\end{align*}

The next result, which is easy to verify, determines the singular vectors
responsible for the composition series of
$M^\mfrak{g}_\mfrak{b}(\lambda)$.
It can be regarded as a degenerate case in the series of parabolic subalgebras
with Heisenberg type nilradicals discussed in \cite{Krizka-Somberg2015}.
We write the statement
for general highest weight $\lambda\in\mathfrak{h}^*$, so that the number of singular
vectors is reduced for the weights which are not dominant and regular.
In particular, in the case $\lambda$ is a regular dominant integral weight there
are six singular vectors.
\medskip

\lemma{Let $\lambda=\lambda_1\omega_1+\lambda_2\omega_2$ be the highest weight for the Verma module $M^\mfrak{g}_\mfrak{b}(\lambda) \simeq \C[\partial_x,\partial_y,\partial_z]$. Then the singular vectors and their weights are
\begin{enumerate}
  \item[$1)$] $v_\lambda=1$, $(\lambda_1, \lambda_2)$
  \item[$2)$] $v_{s_1\cdot \lambda}=\partial_x^{\lambda_1+1}$, $(-\lambda_1-2, \lambda_1+\lambda_2+1)$,
  \item[$3)$] $v_{s_2\cdot \lambda}=\partial_y^{\lambda_2+1}$, $(\lambda_1+\lambda_2+1, -\lambda_2-2)$,
  \item[$4)$] $v_{s_2s_1\cdot \lambda}=\sum_{k=0}^{\lambda_2+1}\!{k! \over 2^k}\binom{\lambda_2+1}{k}
\binom{\lambda_1+\lambda_2+2}{k}\partial_z^k\partial_y^{\lambda_2-k+1}\partial_x^{\lambda_1+\lambda_2-k+2}$,
$(-\lambda_1-\lambda_2-3, \lambda_1)$,
  \item[$5)$] $v_{s_1s_2\cdot \lambda}=\sum_{k=0}^{\lambda_1+1}\!{(-1)^kk! \over 2^k}\binom{\lambda_1+1}{k}
\binom{\lambda_1+\lambda_2+2}{k}\partial_z^k\partial_x^{\lambda_1-k+1}\partial_y^{\lambda_1+\lambda_2-k+2}$,
$(\lambda_2, -\lambda_1-\lambda_2-3)$,
  \item[$6)$] $v_{s_1s_2s_1\cdot \lambda}=\sum_{k=0}^{\lambda_1+\lambda_2+2}\!{k! \over 2^k}\binom{\lambda_1+\lambda_2+2}{k}
\sum_{\ell=0}^k\binom{\lambda_2+1}{\ell}\binom{\lambda_1+1}{k-\ell}(-1)^\ell
\partial_z^k\partial_x^{\lambda_1+\lambda_2-k+2}\partial_y^{\lambda_1+\lambda_2-k+2}$, \\
$(-\lambda_2-2, -\lambda_1-2)$.
\end{enumerate}}

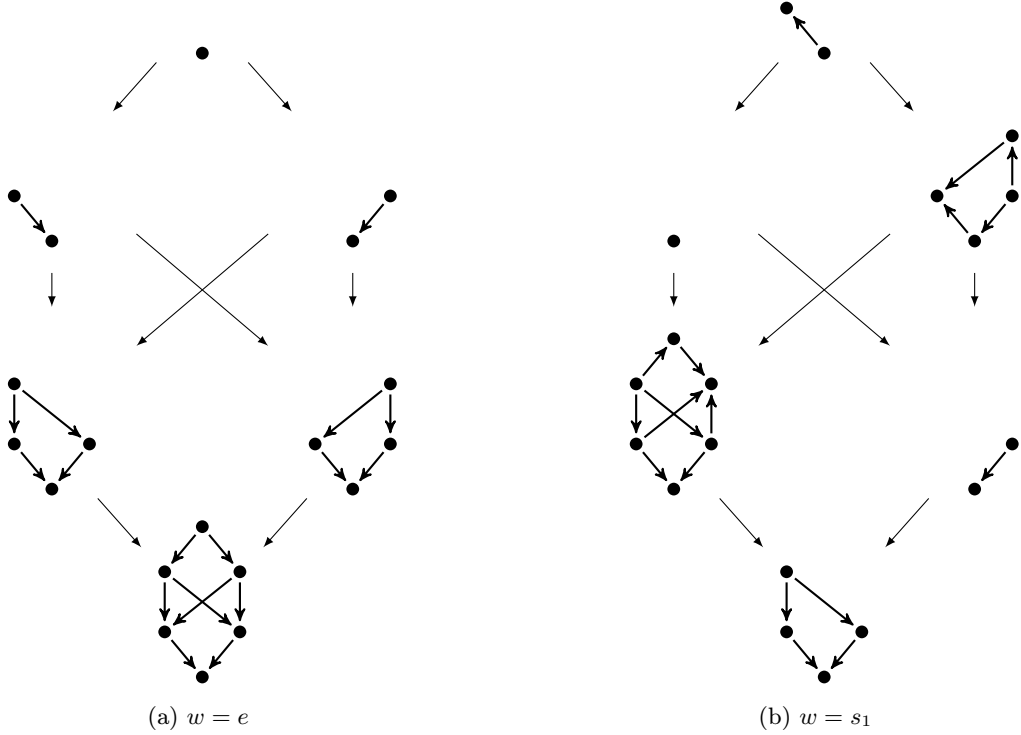
\begin{figure}[t]
\vspace{-0.9cm}
\centering
\subcaptionbox{$w=e$\label{fig:e twisted}}[0.4\textwidth]
{\begin{tikzpicture}
[yscale=1,xscale=1,vector/.style={circle,draw=white,fill=black,ultra thick, inner sep=0.8mm},vector2/.style={circle,draw=white,fill=white,ultra thick, inner sep=1mm}]
\begin{scope}
  \node (A) at (0,0) [vector]  {};
  \node (B) at (0.5,0.6) [vector]  {};
  \node (C) at (-0.5,0.6) [vector]  {};
  \node (D) at (0.5,1.4) [vector]  {};
  \node (E) at (-0.5,1.4) [vector]  {};
  \node (F) at (0,2.0) [vector]  {};
  \node (D0) at (0.7,1.6) {};
  \node (E0) at (-0.7,1.6) {};
  \draw [thick, stealth'-] (A) -- (B);
  \draw [thick, stealth'-] (A) -- (C);
  \draw [thick, stealth'-] (B) -- (D);
  \draw [thick, stealth'-] (C) -- (E);
  \draw [thick, stealth'-] (E) -- (F);
  \draw [thick, stealth'-] (D) -- (F);
  \draw [thick, stealth'-] (C) -- (D);
  \draw [thick, stealth'-] (B) -- (E);
\end{scope}
\begin{scope}[shift={(2,2.5)}]
  \node (A) at (0,0) [vector]  {};
  \node (B) at (0.5,0.6) [vector]  {};
  \node (C) at (-0.5,0.6) [vector]  {};
  \node (D) at (0.5,1.4) [vector]  {};
  \node (E) at (-0.5,1.4) [vector2]  {};
  \node (F) at (0,2.0) [vector2]  {};
  \node (F1) at (0,2.3) {};
  \node (E1) at (-1,1.8) {};
  \node (C1) at (-0.5,0) {};
  \draw [thick, stealth'-] (A) -- (B);
  \draw [thick, stealth'-] (A) -- (C);
  \draw [thick, stealth'-] (B) -- (D);
  \draw [thick, stealth'-] (C) -- (D);
\end{scope}
\begin{scope}[shift={(-2,2.5)}]
  \node (A) at (0,0) [vector]  {};
  \node (B) at (0.5,0.6) [vector]  {};
  \node (C) at (-0.5,0.6) [vector]  {};
  \node (D) at (0.5,1.4) [vector2]  {};
  \node (E) at (-0.5,1.4) [vector]  {};
  \node (F) at (0,2.0) [vector2]  {};
  \node (F2) at (0,2.3) {};
  \node (D2) at (1,1.8) {};
  \node (B2) at (0.5,0) {};
  \draw [thick, stealth'-] (A) -- (B);
  \draw [thick, stealth'-] (A) -- (C);
  \draw [thick, stealth'-] (C) -- (E);
  \draw [thick, stealth'-] (B) -- (E);
\end{scope}
\begin{scope}[shift={(-2,5.8)}]
  \node (A) at (0,0) [vector]  {};
  \node (B) at (0.5,0.6) [vector2]  {};
  \node (C) at (-0.5,0.6) [vector]  {};
  \node (D) at (0.5,1.4) [vector2]  {};
  \node (E) at (-0.5,1.4) [vector2]  {};
  \node (F) at (0,2.0) [vector2]  {};
  \node (A3) at (0,-0.3) {};
  \node (B3) at (1,0.2) {};
  \node (D3) at (0.7,1.6) {};
  \draw [thick, stealth'-] (A) -- (C);
\end{scope}
\begin{scope}[shift={(2,5.8)}]
  \node (A) at (0,0) [vector]  {};
  \node (B) at (0.5,0.6) [vector]  {};
  \node (C) at (-0.5,0.6) [vector2]  {};
  \node (D) at (0.5,1.4) [vector2]  {};
  \node (E) at (-0.5,1.4) [vector2]  {};
  \node (F) at (0,2.0) [vector2]  {};
  \node (A4) at (0,-0.3) {};
  \node (C4) at (-1,0.2) {};
  \node (E4) at (-0.7,1.6) {};
  \draw [thick, stealth'-] (A) -- (B);
\end{scope}
\begin{scope}[shift={(0,8.3)}]
  \node (A) at (0,0) [vector]  {};
  \node (B) at (0.5,0.6) [vector2]  {};
  \node (C) at (-0.5,0.6) [vector2]  {};
  \node (D) at (0.5,1.4) [vector2]  {};
  \node (E) at (-0.5,1.4) [vector2]  {};
  \node (F) at (0,2.0) [vector2]  {};
  \node (C5) at (-0.5,0) {};
  \node (B5) at (0.5,0) {};
\end{scope}
\draw[-latex] (A4) -- (F1);
\draw[-latex] (A3) -- (F2);
\draw[-latex] (C4) -- (D2);
\draw[-latex] (B3) -- (E1);
\draw[-latex] (C1) -- (D0);
\draw[-latex] (B2) -- (E0);
\draw[-latex] (B5) -- (E4);
\draw[-latex] (C5) -- (D3);
\end{tikzpicture}
}
\hfill
\subcaptionbox{$w=s_1$\label{fig:s1 twisted}}[0.5\textwidth]
{
\begin{tikzpicture}
[yscale=1,xscale=1,vector/.style={circle,draw=white,fill=black,ultra thick, inner sep=0.8mm},vector2/.style={circle,draw=white,fill=white,ultra thick, inner sep=1mm}]
\begin{scope}
  \node (A) at (0,0) [vector]  {};
  \node (B) at (0.5,0.6) [vector]  {};
  \node (C) at (-0.5,0.6) [vector]  {};
  \node (D) at (0.5,1.4) [vector2]  {};
  \node (E) at (-0.5,1.4) [vector]  {};
  \node (F) at (0,2.0) [vector2]  {};
  \node (D0) at (0.7,1.6) {};
  \node (E0) at (-0.7,1.6) {};
  \draw [thick, stealth'-] (A) -- (B);
  \draw [thick, stealth'-] (A) -- (C);
  \draw [thick, stealth'-] (C) -- (E);
  \draw [thick, stealth'-] (B) -- (E);
\end{scope}
\begin{scope}[shift={(2,2.5)}]
  \node (A) at (0,0) [vector]  {};
  \node (B) at (0.5,0.6) [vector]  {};
  \node (C) at (-0.5,0.6) [vector2]  {};
  \node (D) at (0.5,1.4) [vector2]  {};
  \node (E) at (-0.5,1.4) [vector2]  {};
  \node (F) at (0,2.0) [vector2]  {};
  \node (F1) at (0,2.3) {};
  \node (E1) at (-1,1.8) {};
  \node (C1) at (-0.5,0) {};
  \draw [thick, stealth'-] (A) -- (B);
\end{scope}
\begin{scope}[shift={(-2,2.5)}]
  \node (A) at (0,0) [vector]  {};
  \node (B) at (0.5,0.6) [vector]  {};
  \node (C) at (-0.5,0.6) [vector]  {};
  \node (D) at (0.5,1.4) [vector]  {};
  \node (E) at (-0.5,1.4) [vector]  {};
  \node (F) at (0,2.0) [vector]  {};
  \node (F2) at (0,2.3) {};
  \node (D2) at (1,1.8) {};
  \node (B2) at (0.5,0) {};
  \draw [thick, stealth'-] (A) -- (B);
  \draw [thick, stealth'-] (A) -- (C);
  \draw [thick, -stealth'] (B) -- (D);
  \draw [thick, stealth'-] (C) -- (E);
  \draw [thick, -stealth'] (E) -- (F);
  \draw [thick, stealth'-] (D) -- (F);
  \draw [thick, -stealth'] (C) -- (D);
  \draw [thick, stealth'-] (B) -- (E);
\end{scope}
\begin{scope}[shift={(-2,5.8)}]
  \node (A) at (0,0) [vector]  {};
  \node (B) at (0.5,0.6) [vector2]  {};
  \node (C) at (-0.5,0.6) [vector2]  {};
  \node (D) at (0.5,1.4) [vector2]  {};
  \node (E) at (-0.5,1.4) [vector2]  {};
  \node (F) at (0,2.0) [vector2]  {};
  \node (A3) at (0,-0.3) {};
  \node (B3) at (1,0.2) {};
  \node (D3) at (0.7,1.6) {};
\end{scope}
\begin{scope}[shift={(2,5.8)}]
  \node (A) at (0,0) [vector]  {};
  \node (B) at (0.5,0.6) [vector]  {};
  \node (C) at (-0.5,0.6) [vector]  {};
  \node (D) at (0.5,1.4) [vector]  {};
  \node (E) at (-0.5,1.4) [vector2]  {};
  \node (F) at (0,2.0) [vector2]  {};
  \node (A4) at (0,-0.3) {};
  \node (C4) at (-1,0.2) {};
  \node (E4) at (-0.7,1.6) {};
  \draw [thick, stealth'-] (A) -- (B);
  \draw [thick, -stealth'] (A) -- (C);
  \draw [thick, -stealth'] (B) -- (D);
  \draw [thick, stealth'-] (C) -- (D);
\end{scope}
\begin{scope}[shift={(0,8.3)}]
  \node (A) at (0,0) [vector]  {};
  \node (B) at (0.5,0.6) [vector2]  {};
  \node (C) at (-0.5,0.6) [vector]  {};
  \node (D) at (0.5,1.4) [vector2]  {};
  \node (E) at (-0.5,1.4) [vector2]  {};
  \node (F) at (0,2.0) [vector2]  {};
  \node (C5) at (-0.5,0) {};
  \node (B5) at (0.5,0) {};
  \draw [thick, -stealth'] (A) -- (C);
\end{scope}
\draw[-latex] (A4) -- (F1);
\draw[-latex] (A3) -- (F2);
\draw[-latex] (C4) -- (D2);
\draw[-latex] (B3) -- (E1);
\draw[-latex] (C1) -- (D0);
\draw[-latex] (B2) -- (E0);
\draw[-latex] (B5) -- (E4);
\draw[-latex] (C5) -- (D3);
\end{tikzpicture}
}
\caption{Generalized weak BGG resolution for $w=e$ and $w=s_1$}
\label{fig:e and s1 twisted}
\end{figure}

The Hasse diagram corresponding to the affine orbit of the Weyl group $W$
for a regular dominant integral weight $\lambda$ of $\mfrak{sl}(3,\C)$ is
drawn on Figure \ref{fig:e twisted}. The nodes of the overall graph correspond to Verma modules
and the arrows are their homomorphisms, and the dots and arrows in each node
(corresponding to a Verma module)
represent the singular vectors and the Verma submodules they generate, respectively.

\subsection{Twisted Verma modules for $w=s_1$}

For $w=s_1$, we have $I_w=(\partial_x, y, z)$, $M^w_\mfrak{g}(\lambda)\simeq \C[x, \partial_y, \partial_z]$, and
\begin{align}
\pi^w_{(\lambda_1,\lambda_2)}=\pi_{(-\lambda_1-1, \lambda_1+\lambda_2+2)}\circ \Ad(\dot{w}^{-1})
\end{align}
since $w^{-1}(\lambda+\rho)=w^{-1}(\lambda_1+1,\lambda_2+1)= (-\lambda_1-1, \lambda_1+\lambda_2+2)$ for $\lambda=\lambda_1\omega_1+\lambda_2\omega_2$.
Because
\begin{align}
\begin{gathered}
\Ad(\dot{w}^{-1})(e_1)=-f_1,\qquad \Ad(\dot{w}^{-1})(e_{12})=e_2, \qquad \Ad(\dot{w}^{-1})(e_2)=-e_{12}, \\
\Ad(\dot{w}^{-1})(f_1)=-e_1, \qquad \Ad(\dot{w}^{-1})(f_{12})=f_2, \qquad
\Ad(\dot{w}^{-1})(f_2)=-f_{12}, \\
\Ad(\dot{w}^{-1})(h_1)=-h_1, \qquad \Ad(\dot{w}^{-1})(h_2)=h_1+h_2, \\
\end{gathered}
\end{align}
we obtain
\begin{align*}
\begin{aligned}
\pi^w_\lambda(f_1)&=-x^2\partial_x - (z-{\textstyle {1\over 2}}xy)\partial_y-({\textstyle {1\over 4}}x^2y+ {\textstyle {1\over 2}}xz)\partial_z + \lambda_1x, \\
\pi^w_\lambda(f_2)&=\partial_z, \\
\pi^w_\lambda(f_{12})&=-\partial_y-{\textstyle {1 \over 2}}x\partial_z, \\
\pi^w_\lambda(e_1)&=\partial_x-{\textstyle {1 \over 2}}y\partial_z, \\
\pi^w_\lambda(e_2)&=-(xz+{\textstyle {1\over 2}}x^2y)\partial_x - (yz-{\textstyle {1\over 2}}xy^2)\partial_y - (z^2+{\textstyle {1 \over 4}}x^2y^2)\partial_z - (\lambda_2+3)z + {\textstyle {1\over 2}}(2\lambda_1+\lambda_2+3)xy, \\
\pi^w_\lambda(e_{12})&=y^2\partial_y -(z+{\textstyle {1\over 2}}xy)\partial_x -({\textstyle {1\over 4}}xy^2-{\textstyle {1\over 2}}yz)\partial_z + (\lambda_1+\lambda_2+3)y, \\
\pi^w_\lambda(h_1)&=-2x\partial_x+y\partial_y-z\partial_z+\lambda_1, \\
\pi^w_\lambda(h_2)&=x\partial_x+y\partial_y+2z\partial_z+\lambda_2+3.
\label{eq:s1 twisted}
\end{aligned}
\end{align*}
The vector $1\in \C[x, \partial_y, \partial_z]$ has the weight $\lambda=(\lambda_1, \lambda_2)$.
In the partial Fourier dual picture of the representation, the Lie algebra $\mfrak{g}$  acts on
the polynomial algebra $\C[x, \xi_y, \xi_z]$ by
\begin{align*}
\begin{aligned}
\hat{\pi}^w_\lambda(f_1)&=\xi_y\partial_{\xi_z}
-x(x\partial_x-{\textstyle {1 \over 2}}\xi_z\partial_{\xi_z}-\lambda_1)
-{\textstyle {1 \over 2}}x(\xi_{y}-{\textstyle {1 \over 2}}x\xi_z)\partial_{\xi_y}, \\
\hat{\pi}^w_\lambda(f_{2})&={\xi_z}, \\
\hat{\pi}^w_\lambda(f_{12})&=-\xi_y-{\textstyle {1 \over 2}}x{\xi_z}, \\
\hat{\pi}^w_\lambda(e_1)&=\partial_x+{\textstyle {1 \over 2}}\xi_z\partial_{\xi_y}, \\
\hat{\pi}^w_\lambda(e_{2})&=(x\partial_x-\xi_y\partial_{\xi_y}-\xi_z\partial_{\xi_z}
+\lambda_2)\partial_{\xi_z}
+{\textstyle {1 \over 2}}x(x\partial_x+\xi_y\partial_{\xi_y}
-{\textstyle {1 \over 2}}x \xi_z\partial_{\xi_y}-2\lambda_1-\lambda_2-1)\partial_{\xi_y}, \\
\hat{\pi}^w_\lambda(e_{12})&=\partial_x\partial_{\xi_z}
+(\xi_y\partial_{\xi_y}+{\textstyle {1 \over 2}}\xi_z\partial_{\xi_z}-\lambda_1-\lambda_2-{\textstyle {1 \over 2}})
\partial_{\xi_y}
+{\textstyle {1 \over 2}}x(\partial_x
-{\textstyle {1 \over 2}}\xi_{z}\partial_{\xi_y})\partial_{\xi_y}, \\
\hat{\pi}^w_\lambda(h_1)&=-2x\partial_x-\xi_y\partial_{\xi_y}+\xi_z\partial_{\xi_z}+\lambda_1, \\
\hat{\pi}^w_\lambda(h_2)&=x\partial_x-\xi_y\partial_{\xi_y}-2\xi_z\partial_{\xi_z}+\lambda_2.
\end{aligned}
\end{align*}

\lemma{\label{lem2.2}The twisted Verma module $\smash{M^w_\mfrak{g}}(\lambda)$ for
$\lambda=\lambda_1\omega_1+\lambda_2\omega_2$ is generated by one vector
$v_\lambda$. For $\lambda_1 \notin \N_0$ this generator is
$v_\lambda=1 \in \C[x,\partial_y,\partial_z]$, while for
$\lambda_1 \in \N_0$ the generator is equal to $v_\lambda=x^{\lambda_1+1} \in \C[x,\partial_y,\partial_z]$.}

\proof{We observe $\pi^w_\lambda(f_1^k)1=k!\binom{\lambda_1}{k}x^k$ for $k\in \N_0$.
It follows that the vectors $\{\pi^w_\lambda(f_1^k)1;\, k \in \N_0\}$ generate the subspace
$\C[x] \subset \C[x,\partial_y,\partial_z]$ for $\lambda_1 \notin \N_0$.
Further, for $\lambda_1 \in \N_0$ we have
$\pi^w_\lambda(e_1^k)x^{\lambda_1+1}=k!\binom{\lambda_1+1}{k}x^{\lambda_1+1-k}$
and $\pi^w_\lambda(f_1^k)x^{\lambda_1+1}=(-1)^kk!x^{\lambda_1+1+k}$ for $k\in \N_0$,
and therefore the vectors
$\{\pi^w_\lambda(f_1^k)x^{\lambda_1+1},\, \pi^w_\lambda(e_1^k)x^{\lambda_1+1};\, k \in \N_0\}$
generate again the subspace $\C[x] \subset \C[x,\partial_y,\partial_z]$. Now, from the form of
elements $\pi^w_\lambda(f_2)=\partial_z$ and $\pi^w_\lambda(f_{12})=-\partial_y-{1 \over 2}x\partial_z$
the rest of the proof easily follows.}

\lemma{\label{lem:generators}Let $\lambda=\lambda_1\omega_1+\lambda_2\omega_2$ and let $I^w_\lambda$
be the left ideal of $U(\mfrak{g})$ defined by
\begin{align}
  I^w_\lambda=\begin{cases}
    (e_1,e_2,e_{12},h_1-\lambda_1,h_2-\lambda_2) & \text{if $\lambda_1 \notin \N_0$}, \\
    (e_2,e_{12},e_1^{\lambda_1+2},f_1e_1,h_1+\lambda_1+2,h_2-\lambda_1-\lambda_2-1) & \text{if $\lambda_1 \in \N_0$}.
  \end{cases}
\end{align}
Then we have $U(\mfrak{g})/I^w_\lambda \simeq M^w_\mfrak{g}(\lambda)$ as $\mfrak{g}$-modules.}

\proof{Let us consider $w_\lambda = 1\ {\rm mod}\ I^w_\lambda$ as an element in $U(\mfrak{g})/I^w_\lambda$.
 The generator $v_\lambda$ of $\smash{M^w_\mfrak{g}}(\lambda)$ constructed in Lemma \ref{lem2.2}
 allows to define a homomorphism $\varphi \colon U(\mfrak{g}) \rarr \smash{M^w_\mfrak{g}}(\lambda)$ of $U(\mfrak{g})$-modules
 by $\varphi(1)=v_\lambda$. Since the generator $v_\lambda$ is annihilated by the left ideal
$I^w_\lambda$ of $U(\mfrak{g})$, we get the surjective homomorphism
$\tilde{\varphi} \colon U(\mfrak{g})/I^w_\lambda \rarr \smash{M^w_\mfrak{g}}(\lambda)$ of $U(\mfrak{g})$-modules.
Then $\tilde{\varphi}$ is an isomorphism once we prove that the modules have the same
formal characters. There are clearly two complementary cases to be considered:
\begin{enum}
\item[i)]
Let us assume first $\lambda_1 \notin \N_0$. Then $h_1w_\lambda=\lambda_1w_\lambda$ and $h_2w_\lambda=\lambda_2w_\lambda$, therefore $U(\mfrak{h})w_\lambda=\C w_\lambda$. By $e_1w_\lambda=0$, $e_2w_\lambda=0$ and $e_{12}w_\lambda=0$ it
follows $U(\mfrak{n})w_\lambda=\C w_\lambda$. The PBW theorem applied to
$\mfrak{g} = \widebar{\mfrak{n}} \oplus \mfrak{h} \oplus \mfrak{n}$ is equivalent to
$U(\mfrak{g})= U(\widebar{\mfrak{n}}) U(\mfrak{n}) U(\mfrak{h})$, hence we get
$U(\mfrak{g})/I^w_\lambda=U(\mfrak{g})w_\lambda= U(\widebar{\mfrak{n}})w_\lambda$.
But the characters of $U(\mfrak{g})/I^w_\lambda$ and $\smash{M^w_\mfrak{g}}(\lambda)$
are equal, because $U(\mfrak{g})/I^w_\lambda$ is free $U(\widebar{\mfrak{n}})$-module
generated by $w_\lambda$ and hence a Verma module with the highest weight $\lambda$.
\item[ii)]
Let us now assume $\lambda_1 \in \N_0$. Then $h_1w_\lambda=(-\lambda_1-2)w_\lambda$ and $h_2w_\lambda=(\lambda_1+\lambda_2+1)w_\lambda$, therefore $U(\mfrak{h})w_\lambda=\C w_\lambda$. By $e_2w_\lambda=0$ and $e_{12}w_\lambda=0$ it
follows that $U(\mfrak{n}_w \cap \mfrak{n})w_\lambda=\C w_\lambda$. We have $\smash{e_1^{\lambda_1+2}}w_\lambda=0$, and so
$U(\widebar{\mfrak{n}}_w \cap \mfrak{n})w_\lambda=\smash{\bigoplus_{k=0}^{\lambda_1+1}} \C e_1^kw_\lambda$.
Finally, the condition $f_1e_1w_\lambda=0$ implies an equality of vector spaces
\begin{align*}
U(\mfrak{n}_w \cap \widebar{\mfrak{n}})U(\widebar{\mfrak{n}}_w \cap \mfrak{n})w_\lambda = {\textstyle \bigoplus_{k=0}^{\lambda_1} \C e_1^{\lambda_1+1-k}w_\lambda \oplus \bigoplus_{k \in \N_0} \!\C f_1^k w_\lambda}.
\end{align*}
By PBW theorem applied to the vector space decomposition
$\mfrak{g} = (\widebar{\mfrak{n}}_w \cap \widebar{\mfrak{n}}) \oplus (\mfrak{n}_w \cap \widebar{\mfrak{n}}) \oplus (\widebar{\mfrak{n}}_w \cap \mfrak{n}) \oplus (\mfrak{n}_w \cap \mfrak{n}) \oplus \mfrak{h}$ we get
\begin{align*}
  U(\mfrak{g})/I^w_\lambda=U(\mfrak{g})w_\lambda=U(\widebar{\mfrak{n}})w_\lambda \oplus {\textstyle \bigoplus_{k=0}^{\lambda_1}}\, U(\mfrak{n}_w \cap \widebar{\mfrak{n}})\, \C e_1^{\lambda_1+1-k}w_\lambda,
\end{align*}
where we used the fact that $U(\mfrak{g})/I^w_\lambda$ is a free $U(\widebar{\mfrak{n}}_w \cap \widebar{\mfrak{n}})$-module. Hence the characters of $U(\mfrak{g})/I^w_\lambda$ and $\smash{M^w_\mfrak{g}}(\lambda)$ coincide.
\end{enum}
Hence the proof is complete.}

The structure of $s_1$-twisted Verma modules $M^w_\mfrak{g}(\lambda)$
for the highest weights $\lambda \in \mfrak{h}^*$ lying on the affine
orbit of $W$ for an integral dominant weight is shown
on the list below. Namely, the points on the figures denote singular
vectors in $M^w_\mfrak{g}(\lambda)$ or some of its quotients, where the vectors
which were not singular vectors in the former twisted Verma module are
singular vectors in a quotient Verma module.
The composition structure among these vectors is expressed by the arrows: there
is a directed arrow from one vector to another vector if and only if for any
choice of vectors projecting under a quotient homomorphism onto singular
vectors is the latter vector generated by the former vector.

The structure of $M^w_\mfrak{g}(\lambda)$ for a specific highest weight $\lambda \in \mfrak{h}^*$ follows immediately from Lemma \ref{lem:generators} and its proof. If $\lambda_1 \notin \N_0$, then $M^w_\mfrak{g}(\lambda)$ is isomorphic to the Verma module $M^\mfrak{g}_\mfrak{b}(\lambda)$. On the other hand, if $\lambda_1 \in \N_0$, then $M^w_\mfrak{g}(\lambda)/N^w_\mfrak{g}(\lambda)$ is isomorphic to the Verma module $M^\mfrak{g}_\mfrak{b}(s_1\cdot \lambda)$, where $N^w_\mfrak{g}(\lambda)$ is the $\mfrak{g}$-submodule generated by $1 \in \C[x,\partial_y,\partial_z]$.
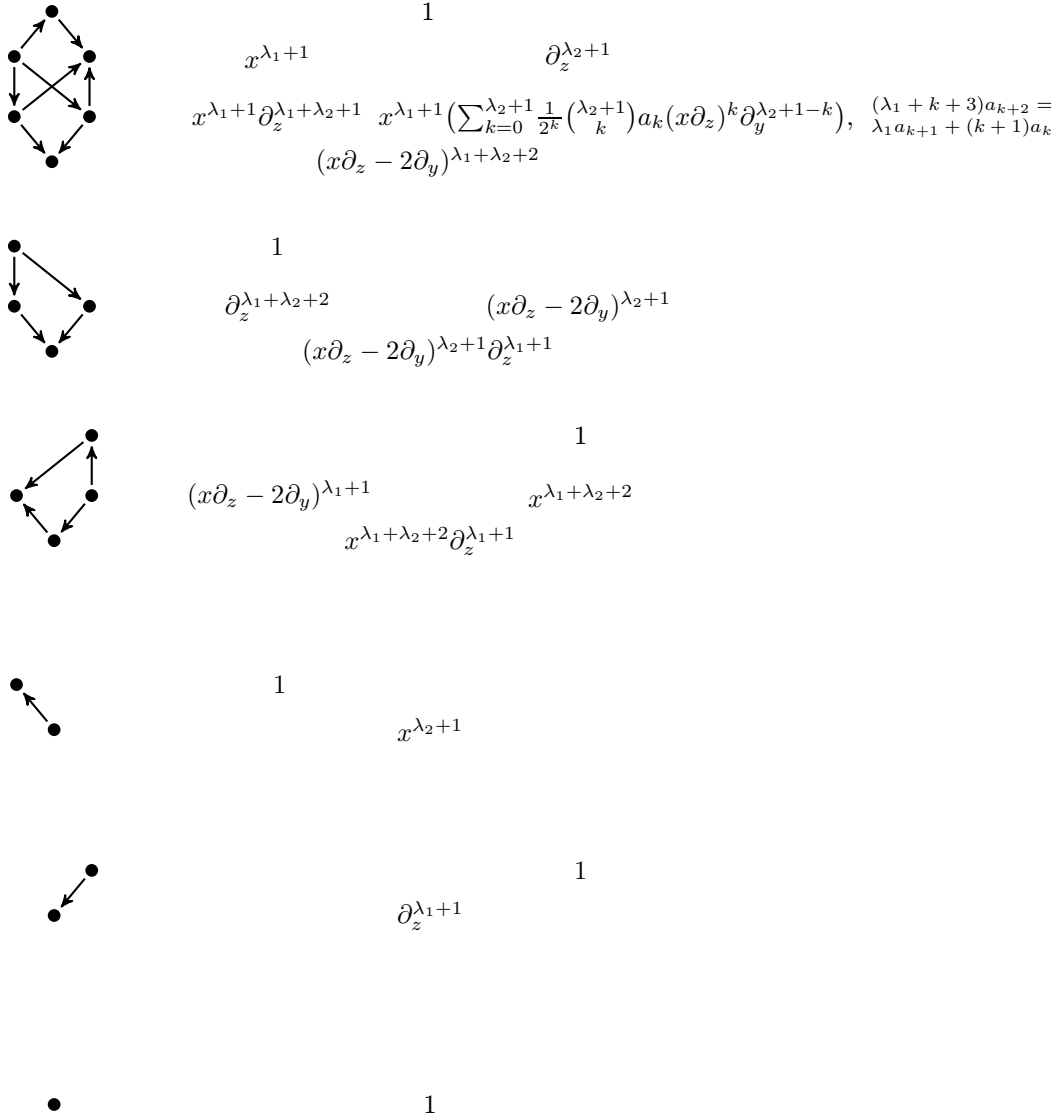
\bigskip

\noindent
\begin{tikzpicture}
[yscale=1,xscale=1,vector/.style={circle,draw=white,fill=black,ultra thick, inner sep=0.8mm},vector2/.style={circle,draw=white,fill=white,ultra thick, inner sep=1mm}]
  \node (A) at (0,0) [vector]  {};
  \node (B) at (0.5,0.6) [vector]  {};
  \node (C) at (-0.5,0.6) [vector]  {};
  \node (D) at (0.5,1.4) [vector]  {};
  \node (E) at (-0.5,1.4) [vector]  {};
  \node (F) at (0,2.0) [vector]  {};
  \draw [thick, stealth'-] (A) -- (B);
  \draw [thick, stealth'-] (A) -- (C);
  \draw [thick, -stealth'] (B) -- (D);
  \draw [thick, stealth'-] (C) -- (E);
  \draw [thick, -stealth'] (E) -- (F);
  \draw [thick, stealth'-] (D) -- (F);
  \draw [thick, -stealth'] (C) -- (D);
  \draw [thick, stealth'-] (B) -- (E);
  \begin{scope}[shift={(5,0)}]
  \node (A) at (0,0) []  {$(x\partial_z-2\partial_y)^{\lambda_1+\lambda_2+2}$};
  \node (B) at (2.5,0.6) []  {$x^{\lambda_1+1} \big(\!\sum_{k=0}^{\lambda_2+1}\! {1 \over 2^k} \binom{\lambda_2+1}{k} a_k(x\partial_z)^k\partial_y^{\lambda_2+1-k}\big)$,};
  \node (B) at (7.1,0.75) [] {\scriptsize{$(\lambda_1+k+3)a_{k+2}=$}};
  \node (B) at (7.1,0.45) [] {\scriptsize{$\lambda_1a_{k+1}+(k+1)a_k$}};
  \node (C) at (-2,0.6) []  {$x^{\lambda_1+1}\partial_z^{\lambda_1+\lambda_2+1}$};
  \node (D) at (2,1.4) []  {$\partial_z^{\lambda_2+1}$};
  \node (E) at (-2,1.4) []  {$x^{\lambda_1+1}$};
  \node (F) at (0,2.0) []  {1};
  \end{scope}
\end{tikzpicture}

\noindent
\begin{tikzpicture}
[yscale=1,xscale=1,vector/.style={circle,draw=white,fill=black,ultra thick, inner sep=0.8mm},vector2/.style={circle,draw=white,fill=white,ultra thick, inner sep=1mm}]
  \node (A) at (0,0) [vector]  {};
  \node (B) at (0.5,0.6) [vector]  {};
  \node (C) at (-0.5,0.6) [vector]  {};
  \node (D) at (0.5,1.4) [vector2]  {};
  \node (E) at (-0.5,1.4) [vector]  {};
  \node (F) at (0,2.0) [vector2]  {};
  \draw [thick, stealth'-] (A) -- (B);
  \draw [thick, stealth'-] (A) -- (C);
  \draw [thick, stealth'-] (C) -- (E);
  \draw [thick, stealth'-] (B) -- (E);
  \begin{scope}[shift={(5,0)}]
  \node (A) at (0,0) []  {$(x\partial_z-2\partial_y)^{\lambda_2+1}\partial_z^{\lambda_1+1}$};
  \node (B) at (2,0.6) []  {$(x\partial_z-2\partial_y)^{\lambda_2+1}$};
  \node (C) at (-2,0.6) []  {$\partial_z^{\lambda_1+\lambda_2+2}$};
  \node (D) at (2,1.4) [vector2]  {};
  \node (E) at (-2,1.4) []  {$1$};
  \node (F) at (0,2.0) [vector2]  {};
  \end{scope}
\end{tikzpicture}

\noindent
\begin{tikzpicture}
[yscale=1,xscale=1,vector/.style={circle,draw=white,fill=black,ultra thick, inner sep=0.8mm},vector2/.style={circle,draw=white,fill=white,ultra thick, inner sep=1mm}]
  \node (A) at (0,0) [vector]  {};
  \node (B) at (0.5,0.6) [vector]  {};
  \node (C) at (-0.5,0.6) [vector]  {};
  \node (D) at (0.5,1.4) [vector]  {};
  \node (E) at (-0.5,1.4) [vector2]  {};
  \node (F) at (0,2.0) [vector2]  {};
  \draw [thick, stealth'-] (A) -- (B);
  \draw [thick, -stealth'] (A) -- (C);
  \draw [thick, -stealth'] (B) -- (D);
  \draw [thick, stealth'-] (C) -- (D);
  \begin{scope}[shift={(5,0)}]
  \node (A) at (0,0) []  {$x^{\lambda_1+\lambda_2+2}\partial_z^{\lambda_1+1}$};
  \node (B) at (2,0.6) []  {$x^{\lambda_1+\lambda_2+2}$};
  \node (C) at (-2,0.6) []  {$(x\partial_z-2\partial_y)^{\lambda_1+1}$};
  \node (D) at (2,1.4) []  {1};
  \node (E) at (-2,1.4) [vector2]  {};
  \node (F) at (0,2.0) [vector2]  {};
  \end{scope}
\end{tikzpicture}

\noindent
\begin{tikzpicture}
[yscale=1,xscale=1,vector/.style={circle,draw=white,fill=black,ultra thick, inner sep=0.8mm},vector2/.style={circle,draw=white,fill=white,ultra thick, inner sep=1mm}]
  \node (A) at (0,0) [vector]  {};
  \node (B) at (0.5,0.6) [vector2]  {};
  \node (C) at (-0.5,0.6) [vector]  {};
  \node (D) at (0.5,1.4) [vector2]  {};
  \node (E) at (-0.5,1.4) [vector2]  {};
  \node (F) at (0,2.0) [vector2]  {};
  \draw [thick, -stealth'] (A) -- (C);
  \begin{scope}[shift={(5,0)}]
  \node (A) at (0,0) []  {$x^{\lambda_2+1}$};
  \node (B) at (2,0.6) [vector2]  {};
  \node (C) at (-2,0.6) []  {$1$};
  \node (D) at (2,1.4) [vector2]  {};
  \node (E) at (-2,1.4) [vector2]  {};
  \node (F) at (0,2.0) [vector2]  {};
  \end{scope}
\end{tikzpicture}

\noindent
\begin{tikzpicture}
[yscale=1,xscale=1,vector/.style={circle,draw=white,fill=black,ultra thick, inner sep=0.8mm},vector2/.style={circle,draw=white,fill=white,ultra thick, inner sep=1mm}]
  \node (A) at (0,0) [vector]  {};
  \node (B) at (0.5,0.6) [vector]  {};
  \node (C) at (-0.5,0.6) [vector2]  {};
  \node (D) at (0.5,1.4) [vector2]  {};
  \node (E) at (-0.5,1.4) [vector2]  {};
  \node (F) at (0,2.0) [vector2]  {};
  \draw [thick, stealth'-] (A) -- (B);
  \begin{scope}[shift={(5,0)}]
  \node (A) at (0,0) []  {$\partial_z^{\lambda_1+1}$};
  \node (B) at (2,0.6) []  {$1$};
  \node (C) at (-2,0.6) [vector2]  {};
  \node (D) at (2,1.4) [vector2]  {};
  \node (E) at (-2,1.4) [vector2]  {};
  \node (F) at (0,2.0) [vector2]  {};
  \end{scope}
\end{tikzpicture}

\noindent
\begin{tikzpicture}
[yscale=1,xscale=1,vector/.style={circle,draw=white,fill=black,ultra thick, inner sep=0.8mm},vector2/.style={circle,draw=white,fill=white,ultra thick, inner sep=1mm}]
  \node (A) at (0,0) [vector]  {};
  \node (B) at (0.5,0.6) [vector2]  {};
  \node (C) at (-0.5,0.6) [vector2]  {};
  \node (D) at (0.5,1.4) [vector2]  {};
  \node (E) at (-0.5,1.4) [vector2]  {};
  \node (F) at (0,2.0) [vector2]  {};
  \begin{scope}[shift={(5,0)}]
  \node (A) at (0,0) []  {1};
  \node (B) at (2,0.6) [vector2]  {};
  \node (C) at (-2,0.6) [vector2]  {};
  \node (D) at (2,1.4) [vector2]  {};
  \node (E) at (-2,1.4) [vector2]  {};
  \node (F) at (0,2.0) [vector2]  {};
  \end{scope}
\end{tikzpicture}
\bigskip

Let us recall that we have the twisting functor $T_w \colon \mcal{O} \rarr \mcal{O}$.
If we apply this functor to the standard BGG resolution for Verma modules shown on
Figure \ref{fig:e twisted}, we obtain twisted BGG resolution shown on Figure \ref{fig:s1 twisted}.

Let us describe the homomorphisms between twisted Verma modules $M^w_\mfrak{g}(\lambda)$ drawn on Figure \ref{fig:s1 twisted} explicitly. Because $\smash{M^w_\mfrak{g}}(\lambda)$ is generated by one vector $v_\lambda \in \smash{M^w_\mfrak{g}}(\lambda)$, a homomorphism $\varphi \colon \smash{M^w_\mfrak{g}}(\lambda) \rarr \smash{M^w_\mfrak{g}}(\mu)$ is uniquely determined by $\varphi(v_\lambda) \in M^w_\mfrak{g}(\mu)$.

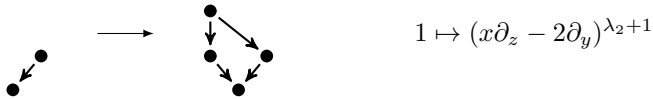
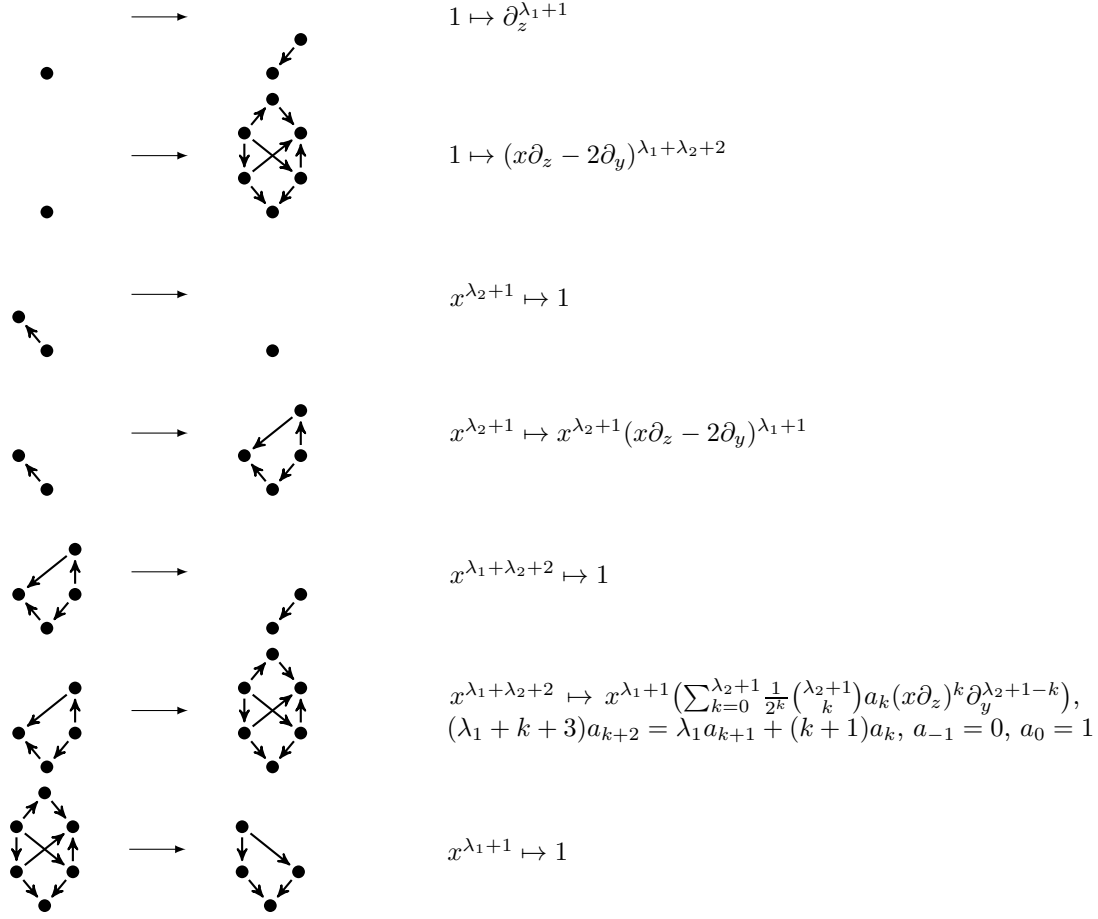
\noindent
\begin{tikzpicture}
[yscale=0.75,xscale=0.75,vector/.style={circle,draw=white,fill=black,ultra thick, inner sep=0.8mm},vector2/.style={circle,draw=white,fill=white,ultra thick, inner sep=1mm}]
\begin{scope}
  \node (A) at (0,0) [vector]  {};
  \node (B) at (0.5,0.6) [vector]  {};
  \node (C) at (-0.5,0.6) [vector2]  {};
  \node (D) at (0.5,1.4) [vector2]  {};
  \node (E) at (-0.5,1.4) [vector2]  {};
  \node (F) at (0,2.0) [vector2]  {};
  \draw [thick, stealth'-] (A) -- (B);
\end{scope}
\begin{scope}[shift={(4,0)}]
  \node (A) at (0,0) [vector]  {};
  \node (B) at (0.5,0.6) [vector]  {};
  \node (C) at (-0.5,0.6) [vector]  {};
  \node (D) at (0.5,1.4) [vector2]  {};
  \node (E) at (-0.5,1.4) [vector]  {};
  \node (F) at (0,2.0) [vector2]  {};
  \draw [thick, stealth'-] (A) -- (B);
  \draw [thick, stealth'-] (A) -- (C);
  \draw [thick, stealth'-] (C) -- (E);
  \draw [thick, stealth'-] (B) -- (E);
\end{scope}
  \draw [-latex] (1.5,1) -- (2.5,1);
  \draw[shift={(7,1)}]
  node [right,text width=0.59\textwidth,inner sep=1mm]
  {$1 \mapsto (x\partial_z-2\partial_y)^{\lambda_2+1}$};
\end{tikzpicture}

\noindent
\begin{tikzpicture}
[yscale=0.75,xscale=0.75,vector/.style={circle,draw=white,fill=black,ultra thick, inner sep=0.8mm},vector2/.style={circle,draw=white,fill=white,ultra thick, inner sep=1mm}]
\begin{scope}
  \node (A) at (0,0) [vector]  {};
  \node (B) at (0.5,0.6) [vector2]  {};
  \node (C) at (-0.5,0.6) [vector2]  {};
  \node (D) at (0.5,1.4) [vector2]  {};
  \node (E) at (-0.5,1.4) [vector2]  {};
  \node (F) at (0,2.0) [vector2]  {};
\end{scope}
\begin{scope}[shift={(4,0)}]
  \node (A) at (0,0) [vector]  {};
  \node (B) at (0.5,0.6) [vector]  {};
  \node (C) at (-0.5,0.6) [vector2]  {};
  \node (D) at (0.5,1.4) [vector2]  {};
  \node (E) at (-0.5,1.4) [vector2]  {};
  \node (F) at (0,2.0) [vector2]  {};
  \draw [thick, stealth'-] (A) -- (B);
\end{scope}
  \draw [-latex] (1.5,1) -- (2.5,1);
  \draw[shift={(7,1)}]
  node [right,text width=0.59\textwidth,inner sep=1mm]
  {$1 \mapsto \partial_z^{\lambda_1+1}$};
\end{tikzpicture}

\noindent
\begin{tikzpicture}
[yscale=0.75,xscale=0.75,vector/.style={circle,draw=white,fill=black,ultra thick, inner sep=0.8mm},vector2/.style={circle,draw=white,fill=white,ultra thick, inner sep=1mm}]
\begin{scope}
  \node (A) at (0,0) [vector]  {};
  \node (B) at (0.5,0.6) [vector2]  {};
  \node (C) at (-0.5,0.6) [vector2]  {};
  \node (D) at (0.5,1.4) [vector2]  {};
  \node (E) at (-0.5,1.4) [vector2]  {};
  \node (F) at (0,2.0) [vector2]  {};
\end{scope}
\begin{scope}[shift={(4,0)}]
  \node (A) at (0,0) [vector]  {};
  \node (B) at (0.5,0.6) [vector]  {};
  \node (C) at (-0.5,0.6) [vector]  {};
  \node (D) at (0.5,1.4) [vector]  {};
  \node (E) at (-0.5,1.4) [vector]  {};
  \node (F) at (0,2.0) [vector]  {};
  \draw [thick, stealth'-] (A) -- (B);
  \draw [thick, stealth'-] (A) -- (C);
  \draw [thick, -stealth'] (B) -- (D);
  \draw [thick, stealth'-] (C) -- (E);
  \draw [thick, -stealth'] (E) -- (F);
  \draw [thick, stealth'-] (D) -- (F);
  \draw [thick, -stealth'] (C) -- (D);
  \draw [thick, stealth'-] (B) -- (E);
\end{scope}
  \draw [-latex] (1.5,1) -- (2.5,1);
  \draw[shift={(7,1)}]
  node [right,text width=0.59\textwidth,inner sep=1mm]
  {$1 \mapsto (x\partial_z-2\partial_y)^{\lambda_1+\lambda_2+2}$};
\end{tikzpicture}

\noindent
\begin{tikzpicture}
[yscale=0.75,xscale=0.75,vector/.style={circle,draw=white,fill=black,ultra thick, inner sep=0.8mm},vector2/.style={circle,draw=white,fill=white,ultra thick, inner sep=1mm}]
\begin{scope}
  \node (A) at (0,0) [vector]  {};
  \node (B) at (0.5,0.6) [vector2]  {};
  \node (C) at (-0.5,0.6) [vector]  {};
  \node (D) at (0.5,1.4) [vector2]  {};
  \node (E) at (-0.5,1.4) [vector2]  {};
  \node (F) at (0,2.0) [vector2]  {};
  \draw [thick, -stealth'] (A) -- (C);
\end{scope}
\begin{scope}[shift={(4,0)}]
  \node (A) at (0,0) [vector]  {};
  \node (B) at (0.5,0.6) [vector2]  {};
  \node (C) at (-0.5,0.6) [vector2]  {};
  \node (D) at (0.5,1.4) [vector2]  {};
  \node (E) at (-0.5,1.4) [vector2]  {};
  \node (F) at (0,2.0) [vector2]  {};
\end{scope}
  \draw [-latex] (1.5,1) -- (2.5,1);
  \draw[shift={(7,1)}]
  node [right,text width=0.59\textwidth,inner sep=1mm]
  {$x^{\lambda_2+1} \mapsto 1$};
\end{tikzpicture}

\noindent
\begin{tikzpicture}
[yscale=0.75,xscale=0.75,vector/.style={circle,draw=white,fill=black,ultra thick, inner sep=0.8mm},vector2/.style={circle,draw=white,fill=white,ultra thick, inner sep=1mm}]
\begin{scope}
  \node (A) at (0,0) [vector]  {};
  \node (B) at (0.5,0.6) [vector2]  {};
  \node (C) at (-0.5,0.6) [vector]  {};
  \node (D) at (0.5,1.4) [vector2]  {};
  \node (E) at (-0.5,1.4) [vector2]  {};
  \node (F) at (0,2.0) [vector2]  {};
  \draw [thick, -stealth'] (A) -- (C);
\end{scope}
\begin{scope}[shift={(4,0)}]
  \node (A) at (0,0) [vector]  {};
  \node (B) at (0.5,0.6) [vector]  {};
  \node (C) at (-0.5,0.6) [vector]  {};
  \node (D) at (0.5,1.4) [vector]  {};
  \node (E) at (-0.5,1.4) [vector2]  {};
  \node (F) at (0,2.0) [vector2]  {};
  \draw [thick, stealth'-] (A) -- (B);
  \draw [thick, -stealth'] (A) -- (C);
  \draw [thick, -stealth'] (B) -- (D);
  \draw [thick, stealth'-] (C) -- (D);
\end{scope}
  \draw [-latex] (1.5,1) -- (2.5,1);
  \draw[shift={(7,1)}]
  node [right,text width=0.59\textwidth,inner sep=1mm]
  {$x^{\lambda_2+1} \mapsto x^{\lambda_2+1}(x\partial_z-2\partial_y)^{\lambda_1+1}$};
\end{tikzpicture}

\noindent
\begin{tikzpicture}
[yscale=0.75,xscale=0.75,vector/.style={circle,draw=white,fill=black,ultra thick, inner sep=0.8mm},vector2/.style={circle,draw=white,fill=white,ultra thick, inner sep=1mm}]
\begin{scope}
  \node (A) at (0,0) [vector]  {};
  \node (B) at (0.5,0.6) [vector]  {};
  \node (C) at (-0.5,0.6) [vector]  {};
  \node (D) at (0.5,1.4) [vector]  {};
  \node (E) at (-0.5,1.4) [vector2]  {};
  \node (F) at (0,2.0) [vector2]  {};
  \draw [thick, stealth'-] (A) -- (B);
  \draw [thick, -stealth'] (A) -- (C);
  \draw [thick, -stealth'] (B) -- (D);
  \draw [thick, stealth'-] (C) -- (D);
\end{scope}
\begin{scope}[shift={(4,0)}]
  \node (A) at (0,0) [vector]  {};
  \node (B) at (0.5,0.6) [vector]  {};
  \node (C) at (-0.5,0.6) [vector2]  {};
  \node (D) at (0.5,1.4) [vector2]  {};
  \node (E) at (-0.5,1.4) [vector2]  {};
  \node (F) at (0,2.0) [vector2]  {};
  \draw [thick, stealth'-] (A) -- (B);
\end{scope}
  \draw [-latex] (1.5,1) -- (2.5,1);
  \draw[shift={(7,1)}]
  node [right,text width=0.59\textwidth,inner sep=1mm]
  {$x^{\lambda_1+\lambda_2+2} \mapsto 1$};
\end{tikzpicture}

\noindent
\begin{tikzpicture}
[yscale=0.75,xscale=0.75,vector/.style={circle,draw=white,fill=black,ultra thick, inner sep=0.8mm},vector2/.style={circle,draw=white,fill=white,ultra thick, inner sep=1mm}]
\begin{scope}
  \node (A) at (0,0) [vector]  {};
  \node (B) at (0.5,0.6) [vector]  {};
  \node (C) at (-0.5,0.6) [vector]  {};
  \node (D) at (0.5,1.4) [vector]  {};
  \node (E) at (-0.5,1.4) [vector2]  {};
  \node (F) at (0,2.0) [vector2]  {};
  \draw [thick, stealth'-] (A) -- (B);
  \draw [thick, -stealth'] (A) -- (C);
  \draw [thick, -stealth'] (B) -- (D);
  \draw [thick, stealth'-] (C) -- (D);
\end{scope}
\begin{scope}[shift={(4,0)}]
  \node (A) at (0,0) [vector]  {};
  \node (B) at (0.5,0.6) [vector]  {};
  \node (C) at (-0.5,0.6) [vector]  {};
  \node (D) at (0.5,1.4) [vector]  {};
  \node (E) at (-0.5,1.4) [vector]  {};
  \node (F) at (0,2.0) [vector]  {};
  \draw [thick, stealth'-] (A) -- (B);
  \draw [thick, stealth'-] (A) -- (C);
  \draw [thick, -stealth'] (B) -- (D);
  \draw [thick, stealth'-] (C) -- (E);
  \draw [thick, -stealth'] (E) -- (F);
  \draw [thick, stealth'-] (D) -- (F);
  \draw [thick, -stealth'] (C) -- (D);
  \draw [thick, stealth'-] (B) -- (E);
\end{scope}
  \draw [-latex] (1.5,1) -- (2.5,1);
  \draw[shift={(7,1)}]
  node [right,text width=0.59\textwidth,inner sep=1mm]
  {$x^{\lambda_1+\lambda_2+2} \mapsto x^{\lambda_1+1} \big(\!\sum_{k=0}^{\lambda_2+1}\! {1 \over 2^k} \binom{\lambda_2+1}{k} a_k(x\partial_z)^k\partial_y^{\lambda_2+1-k}\big)$, $(\lambda_1+k+3)a_{k+2}=\lambda_1 a_{k+1} +(k+1)a_k$, $a_{-1}=0$, $a_0=1$};
\end{tikzpicture}

\noindent
\begin{tikzpicture}
[yscale=0.75,xscale=0.75,vector/.style={circle,draw=white,fill=black,ultra thick, inner sep=0.8mm},vector2/.style={circle,draw=white,fill=white,ultra thick, inner sep=1mm}]
\begin{scope}
  \node (A) at (0,0) [vector]  {};
  \node (B) at (0.5,0.6) [vector]  {};
  \node (C) at (-0.5,0.6) [vector]  {};
  \node (D) at (0.5,1.4) [vector]  {};
  \node (E) at (-0.5,1.4) [vector]  {};
  \node (F) at (0,2.0) [vector]  {};
  \draw [thick, stealth'-] (A) -- (B);
  \draw [thick, stealth'-] (A) -- (C);
  \draw [thick, -stealth'] (B) -- (D);
  \draw [thick, stealth'-] (C) -- (E);
  \draw [thick, -stealth'] (E) -- (F);
  \draw [thick, stealth'-] (D) -- (F);
  \draw [thick, -stealth'] (C) -- (D);
  \draw [thick, stealth'-] (B) -- (E);
\end{scope}
\begin{scope}[shift={(4,0)}]
  \node (A) at (0,0) [vector]  {};
  \node (B) at (0.5,0.6) [vector]  {};
  \node (C) at (-0.5,0.6) [vector]  {};
  \node (D) at (0.5,1.4) [vector2]  {};
  \node (E) at (-0.5,1.4) [vector]  {};
  \node (F) at (0,2.0) [vector2]  {};
  \draw [thick, stealth'-] (A) -- (B);
  \draw [thick, stealth'-] (A) -- (C);
  \draw [thick, stealth'-] (C) -- (E);
  \draw [thick, stealth'-] (B) -- (E);
\end{scope}
  \draw [-latex] (1.5,1) -- (2.5,1);
  \draw[shift={(7,1)}]
  node [right,text width=0.59\textwidth,inner sep=1mm]
  {$x^{\lambda_1+1} \mapsto 1$};
\end{tikzpicture}
\bigskip

There is another aspect of the composition structure of twisted Verma modules, related to a choice of
Lie subalgebras of $\mathfrak{sl}(3,\mathbb{C})$. As for the parallel results for (untwisted) Verma modules,
we refer to \cite{Krizka-Somberg2015}, \cite{koss}. For concreteness,
we shall stick to the case of the Lie algebra $\mathfrak{sl}(2,\mathbb{C})$ embedded on the first
simple root of $\mathfrak{sl}(3,\mathbb{C})$.
\medskip

\lemma{Let us consider the Lie subalgebra $\mathfrak{sl}(2,\mathbb{C})$ of the Lie algebra $\mathfrak{sl}(3,\mathbb{C})$ generated by the elements $\{e_1, f_1, h_1\}$. Then the set of singular vectors for this Lie subalgebra, i.e.\ $v\in M^w_\mfrak{g}(\lambda)$ such that
$\pi^w_\lambda(e_1) v=0$, is given by $\C[\partial_z,x\partial_z-2\partial_y]$, and the corresponding weight spaces with respect to the Cartan subalgebra $\C h_1$ of $\mfrak{sl}(2,\C)$ are
\begin{align}
\bigoplus_{b \in \N_0} \C\partial_z^{a+2b}(x\partial_z-2\partial_y)^b
\end{align}
with the weight $(\lambda_1+a)\omega_1$.}

\proof{First of all, we claim that there is an isomorphism of graded $\C$-algebras
\begin{align}
\C[x,\partial_y,\partial_z]\riso \C[\partial_z,x\partial_z-2\partial_y]\otimes_\C \C[x] \label{isomring}
\end{align}
with $\deg(x)=1$, $\deg(\partial_y)=2$ and $\deg(\partial_z)=1$. The mapping is clearly surjective, and the fact that
$x,\partial_y,\partial_z$ are algebraically independent implies by induction on the degree of $\partial_y$ the
algebraic independence of $x, x\partial_z-2\partial_y, \partial_z$, hence the claim follows.

An elementary calculation shows that any element
in $\C[\partial_z,x\partial_z-2\partial_y]$ is in the kernel of $\pi^w_\lambda(e_1)$,
hence the action of $\pi^w_\lambda(e_1)$ on $\C[\partial_z,x\partial_z-2\partial_y]\otimes_\C \C[x]$ reduces to
the action of $1 \otimes \partial_x$. Therefore, the kernel of $\pi^w_\lambda(e_1)$ on $\C[\partial_z,x\partial_z-2\partial_y]\otimes_\C \C[x]$ is equal to $\C[\partial_z,x\partial_z-2\partial_y]\otimes_\C \C$, and thus to $\C[\partial_z,x\partial_z-2\partial_y]$
by the isomorphism \eqref{isomring}.}

After the application of partial Fourier transform, the isomorphism
\eqref{isomring} can be interpreted as a graded
version of the Fischer tensor product decomposition for the graded algebra
$\mathbb{C}[x,\xi_y,\xi_z]$ with respect to the differential operator
$\partial_x+\frac{1}{2}\xi_z\partial_{\xi_y}$, cf.\ \cite{Fischer1918}. Though we
were not able to find an explicit result on the branching rules for twisted Verma
modules in the available literature, the coincidence of their characters with
characters of (untwisted) Verma modules suggests branching rules in
$K(\mcal{O})$ parallel to those derived in \cite{Kobayashi2012}.

\subsection{Twisted Verma modules for $w=s_1s_2$}

For $w=s_1s_2$, we have $I_w=(x, \partial_y, \partial_z)$, $M^w_\mfrak{g}(\lambda)\simeq \C[\partial_x, y, z]$, and
\begin{align}
\pi^w_{(\lambda_1,\lambda_2)}=\pi_{(\lambda_2+1, -\lambda_1-\lambda_2-2)}\circ \Ad(\dot{w}^{-1})
\end{align}
since $w^{-1}(\lambda+\rho)=w^{-1}(\lambda_1+1, \lambda_2+1)=
(\lambda_2+1, -\lambda_1-\lambda_2-2)$ for $\lambda=\lambda_1\omega_1+\lambda_2\omega_2$.
Due to
\begin{align}
\begin{gathered}
\Ad(\dot{w}^{-1})(e_1)=-f_{12}, \qquad \Ad(\dot{w}^{-1})(e_{12})=-f_2, \qquad \Ad(\dot{w}^{-1})(e_2)=e_1, \\
\Ad(\dot{w}^{-1})(f_1)=-e_{12}, \qquad \Ad(\dot{w}^{-1})(f_{12})=-e_2, \qquad \Ad(\dot{w}^{-1})(f_2)=f_1, \\
\Ad(\dot{w}^{-1})(h_1)=-h_{1}-h_2, \qquad \Ad(\dot{w}^{-1})(h_2)=h_1,
\end{gathered}
\end{align}
we obtain
\begin{align*}
\begin{aligned}
\pi^w_\lambda(f_1)&=-(xz+{\textstyle {1\over 2}}x^2y)\partial_x - (yz-{\textstyle {1\over 2}}xy^2)\partial_y - (z^2+{\textstyle {1 \over 4}}x^2y^2)\partial_z + (\lambda_1-1)z - {\textstyle {1\over 2}}(\lambda_1+2\lambda_2+3)xy, \\
\pi^w_\lambda(f_2)&=-\partial_x+{\textstyle {1 \over 2}}y\partial_z, \\
\pi^w_\lambda(f_{12})&=-y^2\partial_y +(z+{\textstyle {1\over 2}}xy)\partial_x +({\textstyle {1\over 4}}xy^2-{\textstyle {1\over 2}}yz)\partial_z + (\lambda_1+\lambda_2+1)y, \\
\pi^w_\lambda(e_1)&=\partial_z, \\
\pi^w_\lambda(e_2)&=x^2\partial_x + (z-{\textstyle {1\over 2}}xy)\partial_y+({\textstyle {1\over 4}}x^2y+ {\textstyle {1\over 2}}xz)\partial_z + (\lambda_2+2)x, \\
\pi^w_\lambda(e_{12})&=\partial_y+{\textstyle {1 \over 2}}x\partial_z, \\
\pi^w_\lambda(h_1)&=-x\partial_x-y\partial_y-2z\partial_z+\lambda_1-1, \\
\pi^w_\lambda(h_2)&=2x\partial_x-y\partial_y+z\partial_z+\lambda_2+2.
\label{eq:s1s2 twisted}
\end{aligned}
\end{align*}
The vector $1\in \C[\partial_x, y, z]$
has the weight $\lambda=(\lambda_1, \lambda_2)$.
In the partial Fourier dual picture of the representation, $\mathfrak{g}$  acts on
$\C[{\xi_x}, y,z]$ by
\begin{align*}
\begin{aligned}
\hat{\pi}^w_\lambda(f_{1})&=-z(-\xi_x\partial_{\xi_x}+y\partial_y+z\partial_z-\lambda_1)
+{\textstyle {1 \over 2}}y(-\xi_x\partial_{\xi_x}-y\partial_y -{\textstyle {1 \over 2}}y\partial_{\xi_x}\partial_z
+2\lambda_2+\lambda_1+1)\partial_{\xi_x}, \\
\hat{\pi}^w_\lambda(f_2)&=-\xi_x+{\textstyle {1 \over 2}}y\partial_z, \\
\hat{\pi}^w_\lambda(f_{12})&=\xi_x z
-y(y\partial_y+{\textstyle {1 \over 2}}z\partial_z-\lambda_1-\lambda_2-{\textstyle {1 \over 2}})
-{\textstyle {1 \over 2}}y({\xi_{x}}
+{\textstyle {1 \over 2}}y\partial_z)\partial_{\xi_x}, \\
\hat{\pi}^w_\lambda(e_1)&=\partial_z, \\
\hat{\pi}^w_\lambda(e_2)&=z\partial_y
+(\xi_x\partial_{\xi_x}-{\textstyle {1 \over 2}}z\partial_z-\lambda_2)\partial_{\xi_x}
+{\textstyle {1 \over 2}}y(\partial_y+{\textstyle {1 \over 2}}\partial_z\partial_{\xi_x})\partial_{\xi_x}, \\
\hat{\pi}^w_\lambda(e_{12})&=\partial_y-{\textstyle {1 \over 2}}\partial_{\xi_x}\partial_z, \\
\hat{\pi}^w_\lambda(h_1)&=\xi_x\partial_{\xi_x}-y\partial_y-2z\partial_z+\lambda_1, \\
\hat{\pi}^w_\lambda(h_2)&=-2\xi_x\partial_{\xi_x}-y\partial_y+z\partial_z+\lambda_2.
\end{aligned}
\end{align*}

\subsection{Twisted Verma modules for $w=s_1s_2s_1$}

For $w=s_1s_2s_1$, we have $I_w=(\partial_x, \partial_y, \partial_z)$,  $M^w_\mfrak{g}(\lambda)\simeq \C[x, y, z]$, and
\begin{align}
\pi^w_{(\lambda_1,\lambda_2)}=\pi_{(-\lambda_2-1, -\lambda_1-1)}\circ \Ad(\dot{w}^{-1})
\end{align}
since $w^{-1}(\lambda+\rho)=w^{-1}(\lambda_1+1, \lambda_2+1)=
(-\lambda_2-1, -\lambda_1-1)$ for $\lambda=\lambda_1\omega_1+\lambda_2\omega_2$.
Due to
\begin{align}
\begin{gathered}
\Ad(\dot{w}^{-1})(e_1)=-f_2, \qquad \Ad(\dot{w}^{-1})=f_{12}, \qquad \Ad(\dot{w}^{-1})(e_2)=-f_1, \\
\Ad(\dot{w}^{-1})(f_1)=-e_2, \qquad \Ad(\dot{w}^{-1})(f_{12})=e_{12}, \qquad \Ad(\dot{w}^{-1})(f_2)=-e_1, \\
\Ad(\dot{w}^{-1})(h_1)=-h_2, \qquad \Ad(\dot{w}^{-1})(h_2)=-h_1,
\end{gathered}
\end{align}
we obtain
\begin{align*}
\begin{aligned}
\pi^w_\lambda(f_1)&=-y^2\partial_y +(z+{\textstyle {1\over 2}}xy)\partial_x +({\textstyle {1\over 4}}xy^2-{\textstyle {1\over 2}}yz)\partial_z + \lambda_1y, \\
\pi^w_\lambda(f_2)&=-x^2\partial_x - (z-{\textstyle {1\over 2}}xy)\partial_y-({\textstyle {1\over 4}}x^2y+ {\textstyle {1\over 2}}xz)\partial_z + \lambda_2x, \\
\pi^w_\lambda(f_{12})&=(xz+{\textstyle {1\over 2}}x^2y)\partial_x + (yz-{\textstyle {1\over 2}}xy^2)\partial_y + (z^2+{\textstyle {1 \over 4}}x^2y^2)\partial_z - (\lambda_1+\lambda_2)z + {\textstyle {1\over 2}}(\lambda_1-\lambda_2)xy, \\
\pi^w_\lambda(e_1)&=\partial_y+{\textstyle {1 \over 2}}x\partial_z, \\
\pi^w_\lambda(e_2)&=\partial_x-{\textstyle {1 \over 2}}y\partial_z, \\
\pi^w_\lambda(e_{12})&=-\partial_z, \\
\pi^w_\lambda(h_1)&=x\partial_x-2y\partial_y-z\partial_z+\lambda_1, \\
\pi^w_\lambda(h_2)&=-2x\partial_x+y\partial_y-z\partial_z+\lambda_2.
\label{eq:s1s2s1 twisted}
\end{aligned}
\end{align*}
The vector $1\in {\mathbb C}[x, y, z]$
has the weight $\lambda=(\lambda_1, \lambda_2)$.


\section{Outlook and open questions}\label{out}

Let us finish by mentioning that many properties of (untwisted) Verma modules
are not known for twisted Verma modules. For example, it is not clear which
twisted Verma modules are over $U(\mathfrak{g})$ generated by one element. It
is also desirable to understand when a given twisted Verma module
$M^w_{\mathfrak{g}}(\lambda)$ belongs to a parabolic Bernstein-Gelfand-Gelfand
category $\mcal{O}^\mathfrak{p}$ associated to a parabolic subalgebra
$\mathfrak{p}\subset\mathfrak{g}$. This is a non-trivial task, because
the choice of $\mathfrak{p}$ heavily depends on the twisting $w\in W$.
Another question is related to the realization
of twisted Verma modules on Schubert cells -- as for $\mathfrak{sl}(3, \mathbb{C})$ there
are six of them, but there exist altogether eight possibilities for an ideal
defined by annihilating condition for three variables out of the collection
$x,\partial_x,y,\partial_y,z,\partial_z$. What are the isomorphism classes of
the remaining two $\mathfrak{sl}(3, \mathbb{C})$-modules?



\section*{Acknowledgments}

L.\,Křižka is supported by PRVOUK p47,
P.\,Somberg acknowledges the financial support from the grant GA P201/12/G028.



\providecommand{\bysame}{\leavevmode\hbox to3em{\hrulefill}\thinspace}
\providecommand{\MR}{\relax\ifhmode\unskip\space\fi MR }
\providecommand{\MRhref}[2]{%
  \href{http://www.ams.org/mathscinet-getitem?mr=#1}{#2}
}
\providecommand{\href}[2]{#2}

\end{document}